%
%	Invariants of finite group schemes
%	S. Skryabin
%
%	Plain TeX + AMS fonts
%
%	You will need the file  amssym.def  for loading AMS fonts
%	This file is in the standard distribution of amstex
%	and can be downloaded at  ftp.ams.org/pub/tex/amstex/
%
%	PARAMETERS
%
\hsize=5in
\baselineskip=12pt
\vsize=21cm
\parindent=.5cm
\predisplaypenalty=0
\hfuzz=2pt

\def\latexfmt{latex}
\ifx\fmtname\latexfmt\else
%
%	Additional Fonts
%
\input amssym.def
\font\sc=cmcsc10
\font\titlefonts=cmbx12
%
%	Formatting
%
\def\frontmatter{\bgroup
	\leftskip=0pt plus1fil\rightskip=0pt plus1fil
	\parindent=0pt
	\parfillskip=0pt
	\pretolerance=10000}
\def\endfrontmatter{\egroup\bigskip}
\def\title#1{{\baselineskip=1.44\baselineskip\titlefonts#1\par}}
\def\author#1{\bigskip#1\par}
\def\address#1{\bigskip{\it#1}}
\def\thanks#1{\footnote{}{#1\hfil}}
\def\section#1\par{\ifdim\lastskip<\bigskipamount\removelastskip\fi
	\penalty-250\bigskip
	\vbox{\leftskip=0pt plus1fil\rightskip=0pt plus1fil
	\parindent=0pt
	\parfillskip=0pt
	\pretolerance=10000{\bf#1}}\nobreak\medskip}

\fi
\def\proclaim#1. {\medbreak\bgroup{\sc#1.}\hskip.4cm\it}
\def\endproclaim{\egroup
	\ifdim\lastskip<\medskipamount\removelastskip\medskip\fi}
\def\item#1 #2\par{\ifdim\lastskip<\smallskipamount\removelastskip\smallskip\fi
	{\rm#1}\ #2\par\smallskip}
\def\subitem#1 #2\par{\ifdim\lastskip<\smallskipamount\removelastskip\smallskip\fi 
	\indent\phantom{(9)\ }{\rm#1}\ #2\par\smallskip}
\def\Proof#1. {\ifdim\lastskip<\medskipamount\removelastskip\medskip\fi
	{\it Proof.}\quad}
\def\endproof{\quad$\square$\medskip}
\def\Remark. {\ifdim\lastskip<\medskipamount\removelastskip\medskip\fi
	{\it Remark.}\quad}
\def\endremark{\medskip}
%
%	Maths
%
\def\Ad{\mathop{\fam0 Ad}\nolimits}
\def\Aut{\mathop{\fam0 Aut}\nolimits}
\def\can{{\fam0 can}}
\def\chr{\mathop{\fam0 char}\nolimits}
\def\codim{\mathop{\fam0 codim}\nolimits}
\def\degtr{\mathop{\fam0 deg\,tr}\nolimits}
\def\Der{\mathop{\fam0 Der}\nolimits}
\def\End{\mathop{\fam0 End}\nolimits}
\def\Hom{\mathop{\fam0 Hom}\nolimits}
\def\id{\mathop{\fam0 id}\nolimits}
\def\im{\mathop{\fam0 im}}
\def\ind{\mathop{\fam0 ind}\nolimits}
\def\Lie{\mathop{\fam0 Lie}}
\def\rk{\mathop{\fam0 rk}\nolimits}
\def\soc{\mathop{\fam0 soc}\nolimits}
\def\Spec{\mathop{\fam0 Spec}}
\def\iso{\mathrel{\setbox0\hbox{$\rightarrow$}%
  \setbox1\hbox{$\sim$}%
  \dimen0=\wd0\advance\dimen0by-\wd1\advance\dimen0by-.05em%\divide\dimen0by2%
  \box0\kern-\wd1\kern-\dimen0\raise0.75ex\box1\kern\dimen0}}
\def\mapright#1{{}\mathrel{\mathop{\longrightarrow}\limits^{#1}}{}}
\def\longmapright#1#2{{}\mathrel{\smash{\mathop{\count0=#1
  \loop
    \ifnum\count0>0
    \advance\count0 by-1\mathord-\mkern-4mu
  \repeat
  \mathord\rightarrow}\limits^{#2}}}{}}
\def\diagram#1{\vbox{\halign{&\hfil$##$\hfil\cr #1}}}
\def\diagramskip{\noalign{\smallskip}}
\def\trialign#1{\vcenter{\halign{\hfil$\displaystyle##$&\hfil
	$\displaystyle##$\hfil&$\displaystyle##$\hfil\cr #1}}}
\def\smashprod#1{#1\mathbin{\#}\H(G)}
%
%	Symbols
%
\newsymbol\setminus 2272
\newsymbol\nothing 203F
\newsymbol\square 1003
\let\epsilon\varepsilon
\let\phi\varphi
\def\Aop{A^{\fam0op}}
\def\An{{\Bbb A}^n}
\def\a{{\frak a}}
\def\f{{\frak f}}
\def\F{{\cal F}}
\def\g{{\frak g}}
\def\G{{\frak G}}
\def\Gcal{{\cal G}}
\def\h{{\frak h}}
\def\H(#1){k[#1]^*}
\def\M{{\cal M}}
\def\m{{\frak m}}
\def\n{{\frak n}}
\def\O{{\cal O}}
\def\p{{\frak p}}

\def\T{{\cal T}}
\def\Xlr{{X_{\fam0lin.red.}}}
\def\reg{_{\fam0reg}}
\def\greg{{X_{\g{\fam0\mathchar"002Dreg}}}}
\def\Greg{{X\reg}}
\def\Lreg{{X_{L{\fam0\mathchar"002Dreg}}}}
\def\Xreg{{X_{G{\fam0\mathchar"002Dreg}}}}
\def\Mreg{{X_{M{\fam0\mathchar"002Dreg}}}}
\def\Minj{{X_{M{\fam0\mathchar"002Dinj}}}}
%
%		Citations
%
\newcount\citation
\citation=0
\def\citedef#1 {\advance\citation by1
  \expandafter\edef\csname cite:#1\endcsname{{\the\citation}}
  \checkendcitedef}
\def\checkendcitedef#1{\ifx#1\endcitedef\else\citedef#1\fi}
\def\cite#1{[\csname cite:#1\endcsname]}
\citedef
Bou
Cline
DG
Doi
Don
Fri
Hoch
Hoch2
Ja
Kop
Ma
Mil
Mon
Mum
Nich
Ober
Pre
St
Sw
Tak
Veld
Wat
\endcitedef
%
%		Reference section
%
\def\references{\section References\par
	\bgroup
	\parindent=0pt
	%\reffonts
	\rm
	\frenchspacing
	\setbox0\hbox{99. }\leftskip=\wd0
}
\def\endreferences{\egroup}
\newtoks\nextauth
\newif\iffirstauth
\def\checkendauth#1{\ifx\endauth#1%
		\iffirstauth\the\nextauth
		\else{} and \the\nextauth\fi,
	\else\iffirstauth\the\nextauth\firstauthfalse
		\else, \the\nextauth\fi
		\expandafter\auth\expandafter#1\fi}
\def\auth#1 #2 {\nextauth={#1 #2}\checkendauth}
\newif\ifinbook
\newif\ifbookref
\def\nextref#1 {\par\hskip-\leftskip
	\hbox to\leftskip{\hfil\csname cite:#1\endcsname.\ }%
	\initnextref}
\def\initnextref{\bookreffalse\inbookfalse\firstauthtrue\ignorespaces}
\def\paper#1{{#1,}}

\def\book#1{\bookreftrue#1,}
\def\journal#1{#1\ifinbook,\fi}
\def\bookseries#1{#1,}
\def\Vol#1{\ifbookref Vol. #1,\else\ifinbook Vol. #1,\else{\bf#1}\fi\fi}

\def\publisher#1{#1,}
\def\Year#1{\ifbookref #1.\else\ifinbook #1,\else(#1)\fi\fi}
\def\Pages#1{\ifinbook pp. #1.\else #1\expandafter\checktransl\fi}
\long\def\checktransl#1{\ifx[#1\expandafter\takelang\else.#1\fi}
\def\takelang#1]{ [#1].}
%
%		Start of the paper
%
\frontmatter

\title{Invariants of finite group schemes}
\author{Serge Skryabin}
\address{Chebotarev Research Institute, 
%Universitetskaya St.~17, 420008
Kazan, Russia}
%\address{Current address: Max-Planck-Institut f\"ur Mathematik,
%Vivatsgasse 7, 53111 Bonn, Germany}
%\email{Serge.Skryabin@@ksu.ru}
%\email{skryabin@@mpim-bonn.mpg.de}
%\thanks{}

\endfrontmatter

Let $k$ be an algebraically closed field, $G$ a finite group scheme over $k$ 
operating on a scheme $X$ over $k$. Under assumption that $X$ can be covered 
by $G$-invariant affine open subsets the classical results in \cite{DG} 
and \cite{Mum} describe the quotient $X/G$. In case of a free action $X$ is 
known to be a principal homogeneous $G$-space over $X/G$. Furthermore, the 
category of $G$-linearized quasi-coherent sheaves of $\O_X$-modules is 
equivalent then to the category of quasi-coherent sheaves of 
$\O_{X/G}$-modules.

In this paper we attempt to describe the situation when generic stabilizers 
of points on $X$ are nontrivial. To avoid technical complications we assume 
that $X$ is an algebraic variety, although the results 
can be extended to reduced schemes. The stabilizer $G_x$ of a 
rational point $x\in X$ is a subgroup scheme of $G$, and we define its index 
$(G:G_x)$ by analogy with the ordinary finite groups. A point $x$ is {\it 
regular} with respect to the action of $G$ if the index $(G:G_x)$ attains 
the maximal possible value $q(X)$. Theorem 2.1 shows that the set $\Xreg$ of 
all regular points is an open $G$-invariant subset of $X$, the restriction 
to which of the canonical morphism $\pi:X\to X/G$ is finite flat of degree 
$q(X)$. For every $x\in\Xreg$ the fibre $\pi^{-1}\bigl(\pi(x)\bigr)$ is 
$G$-equivariantly isomorphic with the quotient $G_x\backslash G$ and there 
is a bijective correspondence between the $G$-invariant closed subschemes of 
$\Xreg$ and the closed subschemes of $\pi(\Xreg)$. We prove also that the 
field of rational functions $k(X)$ has degree $q(X)$ over the subfield of 
$G$-invariants $k(X)^G$. The arguments used in \cite{DG} and \cite{Mum} are 
essential ingredients in our approach too. At the same time, what we prove 
is not quite a generalization of the classical results as we need more 
restrictions on $X$. 

If $\Xreg=X$ then the equivalence between categories of sheaves mentioned at 
the beginning extends to our settings in only as much as we restrict to 
$G$-linearized quasi-coherent sheaves generated locally by $G$-invariant 
sections (Proposition 3.2). Suppose that $\F$ is an arbitrary $G$-linearized 
coherent sheaf of $\O_X$-modules. In Theorem 3.3 we describe an open 
$G$-invariant subset $U\subset X$ such that the sheaf of $\O_{X/G}$-modules
$(\pi_*\F)^G$ is locally free of rank $s$ over the open subset 
$\pi(U)\subset X/G$, where $s$ is equal to the minimum dimension of the 
subspaces of $G_x$-invariant elements $\F(x)^{G_x}$ in the finite 
dimensional $G_x$-modules $\F(x)=\F_x\otimes_{\O_x}k(x)$, $\,x\in\nobreak U$ 
(here $\F_x$ denotes the stalk of $\F$ at $x$ and $k(x)$ the residue field 
of the local ring $\O_x$). In particular, $(\F\otimes_{\O_X}k(X))^G$ has 
dimension $s$ over $k(X)^G$. We use this result to describe the $G$-socle of 
the $G$-module $\F\otimes_{\O_X}k(X)$ in Corollary 3.4. Given a point 
$x\in\Xreg$ such that $\F_x$ is a free $\O_x$-module, we show in 
Theorem 3.6 that $\F(x)$ is an injective $G_x$-module if and only if there 
exists a $G$-invariant affine open neighbourhood $U$ of $x$ such that 
$\F|_U$ is projective in the category of $G$-linearized sheaves of 
quasi-coherent $\O_U$-modules. Moreover, $\F(U)$ and $\F\otimes_{\O_X}k(X)$ 
are injective $G$-modules in this case. To simplify the statements of 
results we actually consider only affine varieties $X$ and speak about 
modules over the function algebra $k[X]=\O_X(X)$ rather than quasi-coherent 
sheaves. A $G$-linearization on a $k[X]$-module is just a $G$-module 
structure subject to a certain compatibility requirement.

Let us call a group scheme linearly reductive if all its representations are 
completely reducible. Theorem 4.2 says that a point $x\in X$ has a linearly 
reductive stabilizer $G_x$ if and only if $x$ is contained in a 
$G$-invariant affine open subset $U\subset X$ such that $k[U]$ is an 
injective $G$-module. This turns out to be quite a general fact. Unlike 
results in previous sections $X$ can be here any scheme over $k$. When $X$ 
is a variety, the set of points with linearly reductive stabilizers is 
nonempty if and only if $k(X)$ is an injective $G$-module. Moreover, the 
structure of $k(X)$ as a $G$-module is completely determined in this case.

If $\chr k=0$, any finite group scheme over $k$ is constant. Then for all 
$x$ in a nonempty open subset of $X$ the stabilizer $G_x$ coincides with the 
largest subgroup of $G$ acting trivially on the whole $X$. This is the 
reason why our results present an interest mainly for fields of 
characteristic $p>0$. In particular, if $G=\G(\g)$ is the group scheme of 
height one corresponding to a finite dimensional $p$-Lie algebra $\g$ then 
the actions of $G$ on $X$ correspond to actions of $\g$ by derivations of 
the structure sheaf $\O_X$. Probably A.~Milner was the first who observed 
that the degree of $k(X)$ over the subfield of $\g$-invariants $k(X)^\g$ can 
be expressed in terms of Lie algebra stabilizers of points on $X$. He 
considered the special case of the adjoint representation of $\g$ on its 
symmetric algebra $S(\g)$ and used the fact just mentioned to derive a 
lower bound for the maximum dimension of irreducible $\g$-modules 
\cite{Mil}. 

In fact we are able to generalize Theorem 2.1 to the actions of not 
necessarily finite dimensional $p$-Lie algebras (Theorem 5.2). Moreover, if 
$X$ is a smooth affine variety and $f_1,\ldots,f_n$ are $\g$-invariant 
regular functions on $X$, taken in a suitable number, then $k[X]^\g$ is 
generated by $f_1,\ldots,f_n$ over the subalgebra $k[X]^{(p)}$ of $p$-th 
powers in $k[X]$ provided that the differentials $d_xf_1,\ldots,d_xf_n$ 
are linearly independent at all points $x$ in an open 
subset of $X$ whose complement has codimension at least 2 (Theorem 5.4). In 
this case $k[X]^\g$ is free over $k[X]^{(p)}$ and is a locally 
complete intersection ring. A similar result is valid for invariants of 
Frobenius kernels of reduced algebraic groups. This generalizes the work of 
Friedlander and Parshall \cite{Fri}, and Donkin \cite{Don} who considered, 
respectively, the adjoint and the conjugating actions of a semisimple 
algebraic group. We discuss yet another example of the adjoint action of the 
Jacobson-Witt algebra $W_n$. Other applications to invariants of Lie 
algebras of Cartan type will be a subject of separate papers.

I would like to thank the referee for making comments and correction in 
attributing the formula for the $p$-th powers of derivations in section 5.

\section
1. Preliminaries.

Let $k$ be an algebraically closed field. It is the ground field for our 
considerations, so that the functors $\otimes$, $\Hom$ etc.\ are assumed to 
be taken over $k$ unless the base ring is indicated explicitly. Let $G$ be a 
finite group scheme over $k$ and $k[G]$ the associated finite dimensional 
Hopf algebra. We will be considering a group action $\mu:X\times G\to X$ of 
$G$ on a scheme $X$ over $k$ from the right. By \cite{DG} a scheme can be 
regarded as a functor on the category of commutative $k$-algebras. For each 
commutative algebra $K$ the group $G(K)$ operates on $X(K)$, and this action 
is natural in $K$. If $X$ is affine with algebra $k[X]$ then the quotient 
$X/G$ is defined to be $\Spec k[X]^G$ where $k[X]^G\subset k[X]$ is the 
subalgebra of $G$-invariants. More generally, if $X$ can be covered by 
$G$-invariant affine open subsets $U$, then $X/G$ is obtained by patching 
together the affine quotients $U/G$. We list below the properties of the 
canonical morphism $\pi:X\to X/G$ assuming $X$ to be of finite type (see 
\cite{DG}, Ch.~III, \S2, 6.1 and \cite{Mum}, Ch.~III, \S12):

\item
(1) $\pi$ is finite and surjective,

\item
(2) the set-theoretic fibers of $\pi$ coincide with the orbits of the group 
$G(k)$,

\item
(3) $X/G$ has the quotient topology with respect to $\pi$; in particular, 
$\pi$ is open,

\item
(4) if $U\subset X$ is a $G$-invariant open subscheme, then $U/G\cong\pi(U)$ 
and $U=\pi^{-1}\bigl(\pi(U)\bigr)$.

According to \cite{DG}, Ch.~III, \S2, 2.3 the action is said to be free if 
$G(K)$ operates freely on $X(K)$ for each commutative algebra $K$. In case 
of a free action $\pi$ is finite flat of degree $|G|=\dim k[G]$ (which means 
that $k[U]$ are free of rank $|G|$ over $k[U]^G$ for a suitable covering of 
$X$ by $G$-invariant affine open subsets) and the canonical morphism
$\nu=(p_1,\mu):X\times G\to X\times_{X/G}X$,
where $p_1:X\times G\to X$ denotes the projection, is an isomorphism.

If $G'\subset G$ is a subgroup scheme, we let $G'\backslash G$ denote the 
quotient with respect to the action of $G'$ on $G$ by left translations. 
Then $G'\backslash G$ is a finite scheme with the algebra $k[G'\backslash 
G]\cong k[G]^{G'}$, the invariants with respect to the left regular 
representation of $G'$ on $k[G]$. Define the {\it index} of $G'$ in $G$ to 
be $(G:G')=\dim k[G'\backslash G]$. The algebra $k[G]$ is a free module of 
rank $|G'|$ over $k[G'\backslash G]$. Hence $|G|=(G:G')\cdot |G'|$. The 
index can be interpreted in terms of dual Hopf algebras. By \cite{Ober} or 
\cite{Nich} $k[G]^*$ is free both as a left and a right module over its 
subalgebra $k[G']^*$. Clearly, the ranks of these modules are equal to 
$(G:G')$.

If $x\in X(k)$ then its stabilizer $G_x\subset G$ is a subgroup scheme such 
that $G_x(K)=\{g\in G(K)\mid x_Kg=x_K\}$ for each commutative algebra $K$, 
where $x_K$ denotes the image of $x$ in $X(K)$. Let $i_x:\Spec k\to X$ be 
the morphism corresponding to $x$ and
$$
\mu_x:G\cong\Spec k\times G\longmapright5{i_x\times\id_G}X\times 
G\mapright{\mu}X.
$$
the orbit morphism. Then $G_x$ coincides with the fiber of $\mu_x$ above $x$ 
and $\mu_x$ factors through a morphism $\rho:G_x\backslash G\to X$. By 
\cite{DG}, Ch.~III, \S3, 5.2 $\rho$ is an immersion. In fact $\rho$ is a 
closed immersion because $G'\backslash G$ has only rational points. Since 
the case of finite group schemes is especially easy, below we sketch a proof 
of an equivalent assertion for the reader's convenience:

\proclaim
Proposition 1.1.
Suppose that $A\subset k[G]$ is a subalgebra, stable under the right 
regular representation of $G$, and $\m_A=\m\cap A$ where $\m$ is the 
augmentation ideal of $k[G]$. Then $A=k[G'\backslash G]$ where $G'$ is the 
stabilizer of $\m_A$ in $G$.
\endproclaim

\Proof. Put $B=k[G'\backslash G]$, $R=k[G]$. The group scheme $G$ operates 
from the right on $X=\Spec A$ and the inclusion $A\subset R$ corresponds to 
a $G$-equivariant morphism $G\to X$. The latter is the orbit morphism 
$\mu_x$ of the point $x\in X(k)$ corresponding to $\m_A$. Since $G'=G_x$, 
the orbit morphism factors through $G'\backslash G$, which means that 
$A\subset B$. Let $K=R/R\m_A$, and let $g\in G(K)$ be the point 
corresponding to the canonical homomorphism $R\to K$. Since the composite 
homomorphism $A\to R\to K$ factors through $A/\m_A$, we have 
$x_Kg=\mu_x(g)=x_K$, which yields $g\in G'(K)$. It follows that the 
composite $B\to R\to K$ factors through $B/\m_B$ where $m_B=\m\cap B$. In 
other words, $\m_B=R\m_A\cap B$.  However, $R\m_A\cap B=B\m_A$ because $R$ 
is free over $B$. Hence $B=k+\m_B=A+B\m_A$. Since $R$ is finite over $A$, 
the map $G(k)\to X(k)$ determined by $\mu_x$ is surjective. This means that 
$G(k)$ transitively permutes the maximal ideals of $A$. Then $B=A+B\n$ for 
all maximal ideals $\n$ of $A$. An application of Nakayama's lemma yields 
$B=A$.  
\endproof

Suppose now that $A$ is any unital associative algebra on which $G$ operates 
by automorphisms. This means that $A$ has a $G$-module structure and for 
each commutative algebra $K$ the group $G(K)$ operates on $A\otimes K$ via 
automorphisms. We call $A$ a {\it$G$-algebra} in this case. By an 
{\it$(A,G)$-module} we mean a right $A$-module $M$ equipped with an 
additional $G$-module structure such that the $A$-module structure map 
$M\otimes A\to M$ is $G$-equivariant. Denote by $\M_A$ the category of 
$(A,G)$-modules.  The morphisms in $\M_A$ are maps which are simultaneously 
$A$-module and $G$-module homomorphisms.  

This definition is meaningful for an arbitrary group scheme. When $G$ 
is finite, the category of $G$-modules is equivalent to the category of left 
$\H(G)$-modules (see \cite{Ja}, Part I, Ch.~8), and $\M_A$ is equivalent to 
the category of left modules over the smash product algebra $\smashprod\Aop$ 
where $\Aop$ is the algebra $A$ with the opposite multiplication. When $A$ 
is commutative, $\Aop=A$. We refer the reader to \cite{Sw} or \cite{Mon} 
concerning the precise definition of smash products.  Here we just point out 
that $\smashprod{A}$ contains $A$ and $\H(G)$ as subalgebras and the 
multiplication map $A\otimes\H(G)\to\smashprod{A}$ is bijective. One of our 
tools is the following theorem whose interpretations can be found in 
\cite{Cline}, \cite{Kop}, \cite{Tak}:

\proclaim
Imprimitivity Theorem.
Suppose that $A=k[G'\backslash G]$, and let $\m_A$ be the augmentation ideal 
of $A$. Then the functor $M\mapsto M/M\m_A$ is an equivalence between $\M_A$ 
and the category of $G'$-modules. If $M\in\M_A$, then there is an 
isomorphism of $G$-modules $M\cong\ind_{G'}^G M/M\m_A$.  
\endproclaim

Given a $G'$-module $V$, the induced $G$-module is $\ind_{G'}^G V=(V\otimes 
k[G])^{G'}$. To be precise, we compute the $G'$-invariants with respect to 
the tensor product of the $G'$-module structure on $V$ and the left regular 
$G'$-module structure on $k[G]$, and we use the right regular $G$-module 
structure on $k[G]$ to get one on $\ind_{G'}^G V$. This differs from the 
conventions adopted in \cite{Ja}. In terms of the dual algebras 
$$
\ind_{G'}^G V\cong\Hom_{\H(G')}\bigl(\H(G),\,V\bigr).
$$
It follows that $\dim\ind_{G'}^G V=(G:G')\cdot\dim V.$
We mention also that the injective $G$-modules are projective, and vice 
versa, because finite dimensional Hopf algebras are Frobenius (see \cite{Sw}).
Waterhouse proved that every finite group scheme is geometrically reductive 
\cite{Wat}.

Below we prove two lemmas.
Suppose that $A$ is a commutative integral domain and $K$ its field of 
fractions. For a prime ideal $\p$ of $A$ denote $k(\p)=A_\p/A_\p\p$.

\proclaim
Lemma 1.2.
Let $F$ be a finitely generated projective A-module, $F'$ its submodule.  
Denote by $I(\p)$ the image of the canonical map 
$F'\otimes_Ak(\p)\to F\otimes_Ak(\p)$ and put $q=\dim_KF'\otimes_AK$. Then:

\item
(1) $q=\max_{\p\in\Spec A}\dim_{k(\p)}I(\p)$.

\item
(2) The subset $U=\{\p\in\Spec A\mid\dim_{k(\p)}I(\p)=q\}$ is Zariski open 
and coincides with the subset $U'=\{\p\in\Spec A\mid F'_\p$ is a direct 
summand of $F_\p\}$.

\item
(3) If $U=\Spec A$ then $F'$ is projective of rank $q$ and is a direct 
summand of $F$.

\endproclaim

\Proof.
We identify all localizations of $F$, as well as localizations of $F'$, with 
their images in $F\otimes_AK$. Let $\p\in\Spec A$. Suppose that 
$u_1,\ldots,u_r\in F'_\p$ are elements whose images in $F_\p/\p F_\p\cong 
F\otimes_Ak(\p)$ are linearly independent over $k(\p)$. Then $S=A_\p 
u_1+\ldots+A_\p u_r$ is a direct summand of the free $A_\p$-module $F_\p$ 
(see \cite{Bou}, Ch.~II, \S3, Cor.~1 to Prop.~5). In particular, $S$ is 
free over $A_\p$ with a basis $u_1,\ldots,u_r$. We get $r=\dim_K 
S\otimes_AK\le q$. It follows $\dim_{k(\p)}I(\p)\le q$. By definition of $q$ 
the equality holds here for $\p=(0)$.

Suppose that $\p\in U$. Then we can take $q$ elements $u_1,\ldots,u_q$ 
above, so that $r=q$. By modularity law $S$ is a direct summand of 
$F'_\p$. So $F'_\p=S\oplus C$ where $C\subset F_\p$ is an 
$A_\p$-submodule. Tensoring with $K$ and comparing dimensions over $K$, we 
deduce that $C\otimes_AK=0$. However, $F_\p$ is torsion-free since 
$A_\p$ is a domain. It follows $C=0$, i.e., $F'_\p=S$. Thus $\p\in U'$.

Conversely, suppose that $\p\in U'$. Then $F'_\p$ is free of rank $q$ over 
$A_\p$ and the map $F'\otimes_Ak(\p)\to F\otimes_Ak(\p)$ is 
injective. It follows $\dim_{k(\p)}I(\p)=q$. There exists an $A_\p$-module 
epimorphism $\phi:F_\p\to F'_\p$ which is identity on $F'_\p$. Since $F$ is 
finitely generated, $\phi(F)\subset F'_s$ for a suitable $s\in 
A\setminus\p$. Then $\phi(F_s)=F'_s$, and so $F'_s$ is a direct summand of 
$F_s$. This shows that $U'$ is open in $\Spec A$. Finally, (3) is obtained 
by an application of \cite{Bou}, Ch. II, \S5, Th.~1 and \S3, Cor.~1 to Prop. 
12.  
\endproof

\proclaim
Lemma 1.3.
Let $F$ be a finitely generated projective $A$-module, $F'$ and $F''$ its 
direct summands. Denote by $I'(\p)$, $I''(\p)$ the images, respectively, of 
$F'\otimes_Ak(\p)$, $F''\otimes_Ak(\p)$ in $F\otimes_Ak(\p)$, and put 
$s=\min_{\p\in\Spec A}\dim_{k(\p)}I'(\p)\cap I''(\p).$ Then:

\item
(1) The subset $U=\{\p\in\Spec A\mid\dim_{k(\p)}I'(\p)\cap I''(\p)=s\}$ is open 
in $\Spec A$ and consists precisely of those $\p$ for which $F'_\p+F''_\p$ 
is a direct summand of $F_\p$.

\item
(2) For all $\p\in U$ the canonical map $(F'\cap F'')\otimes_Ak(\p)\to 
F\otimes_Ak(\p)$ is an isomorphism onto $I'(\p)\cap I''(\p)$. 

\item
(3) If $U=\Spec A$ then the $A$-module $F'\cap F''$ is projective of rank 
$s$ and is a direct summand of $F$.

\endproclaim

\Proof.
Since $F'$, $F''$ are direct summands of $F$, the dimensions of vector 
spaces $I'(\p)$, $I''(\p)$ do not depend on $\p$. Then $U$ is the set of 
those $\p\in\Spec A$ for which $I'(\p)+I''(\p)$ has maximal possible 
dimension. Apply now Lemma 1.2 taking $F'+F''$ instead of $F'$ in it. We 
get assertion (1) of Lemma 1.3. If $\p\in U$ then $F'_\p+F''_\p$ is a free 
$A_\p$-module.  Since $F''_\p$ is a direct summand of $F'_\p+F''_\p$, the 
$A_\p$-module 
$$ 
F'_\p/(F'_\p\cap F''_\p)\cong(F'_\p+F''_\p)/F''_\p\eqno(*) 
$$ 
is free too. Then $F'_\p\cap F''_\p$ is a direct summand of $F'_\p$, 
hence also a direct summand of $F_\p$. In particular, 
$(F'\cap F'')_\p\cong F'_\p\cap F''_\p$ is free over $A_\p$ and the map in 
(2) is injective. Denoting by $I(\p)$ the image of that map and tensoring 
$(*)$ with $k(\p)$, we obtain an isomorphism 
$I'(\p)/I(\p)\cong\bigl(I'(\p)+I''(\p)\bigr)/I''(\p)$. It follows 
$I(\p)=I'(\p)\cap I''(\p)$. The last assertion is a special case of Lemma 
1.2(3).  
\endproof

\section
2. The set of $G$-regular points and the properties of the quotient.

Let $G$ be a finite group scheme operating from the right on an irreducible 
algebraic variety $X$. Suppose that $X$ can be covered by $G$-invariant 
affine open subsets, so that $X/G$ exists (as is well known it suffices to 
require that the $G(k)$-orbit of each closed point of $X$ is contained in 
an affine open subset). We will be considering only closed points of $X$, 
so that $x\in X$ means $x\in X(k)$. If $U\subset X$ is an open subset, 
stable under all automorphisms of $X$ determined by the elements of $G(k)$, 
then the composite morphism $U\times G\to X\times G\mapright{\mu}X$ factors 
through $U$, i.e., $U$ is $G$-invariant. In particular, if $U\subset X$ is 
any open subset, $\cap_{g\in G(k)}Ug$ is a $G$-invariant open subset 
contained in $U$. It follows that the field of rational functions $k(X)$ is 
a direct limit of the $G$-algebras $k[U]$ where $U$ runs through the 
$G$-invariant affine open subvarieties of $X$. Hence $G$ operates on $k(X)$ 
by automorphisms. Put 
$$ 
q(X)=\max_{x\in X}\,\,(G:G_x),  
$$\removelastskip 
$$
\Xreg=\{x\in X\mid(G:G_x)=q(X)\}.
$$

\proclaim
Theorem 2.1. $\Xreg$ is a $G$-invariant open subset of $X$. Furthermore:

\item
(1) $\pi\vert_\Xreg:\Xreg\to\pi(\Xreg)$ is a finite flat morphism of degree 
$q(X)$.

\item
(2) For every $x\in\Xreg$ the fibre of $\pi$ above $\pi(x)$ is 
$G$-equivariantly isomorphic with $G_x\backslash G$.

\item
(3) The $G$-invariant closed subschemes $Z$ of $\Xreg$ are in a bijective 
correspondence with the closed subschemes $W$ of $\pi(\Xreg)$. If $Z$ and 
$W$ correspond to each other, then $W\cong Z/G$ and $Z\cong W\times_{X/G}X$.

\item
(4) One has $(Y\times_{X/G}X)/G\cong Y$ for any scheme $Y$ on which $G$ 
operates trivially and a morphism $Y\to\pi(\Xreg)$.

\item
(5) $[k(X):k(X)^G]=q(X)$.

\endproclaim

We first reformulate the assertions of the theorem in the affine case. 

\proclaim
Proposition 2.2.
Suppose that $X\cong\Spec A$ and $\Xreg=X$. Then:

\item
(1) $A$ is a projective $A^G$-module of rank $q(X)$.

\item
(2) If $\m$ is a maximal ideal of $A$ and $\n=\m\cap A^G$ then $\n A$ is a 
maximal $G$-invariant ideal of $A$ and the algebra $A/\n A$ is 
$G$-equivariantly isomorphic with $k[G_\m\backslash G]$ where $G_\m$ is the
stabilizer of $\m$ in $G$.

\item
(3) The assignment $I\mapsto I^G$ establishes a bijection between the 
$G$-invariant ideals of $A$ and the ideals of $A^G$. The inverse 
correspondence is given by $J\mapsto JA$. The canonical maps $A^G\to(A/I)^G$ 
are surjective.

\item
(4) If $B$ is an $A^G$-algebra on which $G$ operates trivially then 
$(B\otimes_{A^G}A)^G\cong B$.

\endproclaim

\Proof.
Given any covering of $X$ by $G$-invariant open subvarieties, it suffices to 
prove the theorem for the induced action of $G$ on each of these 
subvarieties. In particular, we may assume $X$ to be affine. Let $A=k[X]$, 
$R=k[G]$, and let $\mu^*:A\to A\otimes R$, $i_x^*:A\to k$, $\mu_x^*:A\to 
R$ be the comorphisms of $\mu,i_x,\mu_x$, respectively. Denote by $\m_x$ the 
maximal ideal of $A$ consisting of functions vanishing at $x$.

Consider $F=A\otimes R$ as an $A$-module by means of the algebra homomorphism 
$p_1^*:A\to A\otimes R$, $\,a\mapsto a\otimes1$. Clearly $F$ is free of 
finite rank over $A$. Put 
$$F'=(A\otimes1)\cdot\mu^*(A)\subset F.$$ 
Then $F'$ is an $A$-submodule of $F$ generated by $\mu^*(A)$. Denote by 
$I(x)$ the image of the canonical map $F'/\m_xF'\to F/\m_xF$. We have 
$F/\m_xF\cong A/\m_x\otimes R\cong R$ and $I(x)=\mu_x^*(A)$ since 
$\mu_x^*=(i_x^*\otimes\id_R)\circ\mu^*$. Now $\mu_x^*$ is a $G$-equivariant 
algebra homomorphism. Since $G_x$ coincides with the stabilizer of $\m_x$ in 
$G$, Proposition 1.1 ensures $\mu_x^*(A)=k[G_x\backslash G]$. Hence $\dim 
I(x)=(G:G_x)$. It follows that $\Xreg$ coincides with the set $U$ of those 
points $x\in X$ for which $\dim I(x)$ attains its maximal value $q=q(X)$. 
We can now apply Lemma 1.2. By (2) of the lemma $U$ is open. Each $g\in 
G(k)$ determines an inner automorphism of $G$ which induces an isomorphism 
$G_x\cong G_{xg}$. Hence $(G:G_{xg})=(G:G_x)$. It follows that $U$ is 
$G$-invariant.

Let $y\in U$ and let $O\subset X$ be the $G(k)$-orbit of $y$. Then $\dim 
I(z)=q$ for all $z\in O$. Since $O$ is finite, we can find 
$a_1,\ldots,a_q\in A$ such that $\mu_z^*(a_1),\ldots,\mu_z^*(a_q)$ are a 
basis of $I(z)$ for each $z\in O$. Furthermore, we may assume $a_1=1$ since 
$\mu_x^*(1)=1$ for all $x$. Applying Lemma 1.2 to the $A$-submodule of 
$F$ generated by $\mu^*(a_1),\ldots,\mu^*(a_q)$, we see that the set $U_1$ 
of those $x\in X$ for which $\mu_x^*(a_1),\ldots,\mu_x^*(a_q)$ are linearly 
independent is open in $X$. Clearly $U_1\subset U$. Since $\pi$ is 
a finite morphism, the set $W=\pi(X\setminus U_1)$ is closed in $X/G$. Since 
$\pi^{-1}\bigl(\pi(y)\bigr)=O\subset U_1$, we have $\pi(y)\notin W$. Let $V$ 
be an open affine neighbourhood of $\pi(y)$ in $X/G$. Then $\pi^{-1}(V)$ is 
an open affine $G$-invariant neighbourhood of $y$ in $X$ and 
$\pi^{-1}(V)\subset U_1$.

To prove the remainder of the theorem we can again use the local character 
of the assertions and pass to the actions of $G$ on the invariant open 
subsets of the form $\pi^{-1}(V)$ constructed above. We may thus assume that 
$U_1=X$.

By Lemma 1.2 $F'$ is a direct summand of the $A$-module $F$ and for each 
maximal ideal $\m$ of $A$ the localization $F'_\m$ is free of rank $q$ over 
$A_\m$. If $\phi:A^q\to F'$ is the $A$-module homomorphism sending the 
standard generators of $A^q$ to $\mu^*(a_1),\ldots,\mu^*(a_q)$ then the 
localizations of $\phi$ at maximal ideals of $A$ are all isomorphisms. Hence 
$\phi$ is itself an isomorphism, i.e., $F'$ is a free $A$-module with a 
basis $\mu^*(a_1),\ldots,\mu^*(a_q)$. Hence for each $a\in A$ there are 
$b_1,\ldots,b_q\in A$ such that
$$
\mu^*(a)=\sum\,(b_i\otimes1)\cdot\mu^*(a_i).\eqno(*)
$$
Let $\epsilon:R\to k$ be the counit and $m^*:R\to R\otimes R$ the 
comultiplication maps. Applying $\id_A\otimes\,\epsilon$ to both sides of 
$(*)$, we get $a=\sum b_ia_i$ since 
$(\id_A\otimes\,\epsilon)\circ\mu^*=\id_A$.  
Applying $\mu^*\otimes\id_R$ and $\id_A\otimes m^*$ to both sides of 
$(*)$, and taking into account the identity 
$(\mu^*\otimes\id_R)\circ\mu^*=(\id_A\otimes m^*)\circ\mu^*$, we get
$$ 
\sum\,(\mu^*(b_i)\otimes1)\cdot(\mu^*\otimes\id_R)\mu^*(a_i)= 
\sum\,(b_i\otimes1\otimes1)\cdot(\mu^*\otimes\id_R)\mu^*(a_i)\eqno(**)
$$
in $A\otimes R\otimes R$. If $\gamma:A\to A'$ is a homomorphism of 
commutative algebras then $A'\otimes_AF'$ is a free $A'$-module with basis 
elements $1\otimes\mu^*(a_i)$, $\,i=1,\ldots,q$. Since the canonical map 
$A'\otimes_AF'\to A'\otimes_AF\cong A'\otimes R$ is injective, the elements 
$$
(\gamma\otimes\id_R)\mu^*(a_i)\in A'\otimes R,\qquad 
i=1,\ldots,q,
$$ 
are linearly independent over $A'$. Taking $A'=A\otimes R$ and $\gamma=\mu^*$, 
we deduce from $(**)$ that $\mu^*(b_i)=b_i\otimes1$, that is, $b_i\in A^G$ 
for all $i$. Hence $A=A^Ga_1+\ldots+A^Ga_q$. If now $\sum c_ia_i=0$ for some 
$c_1,\ldots,c_q\in A^G$ then $\sum\,(c_i\otimes1)\mu^*(a_i)=\mu^*(\sum 
c_ia_i)=0$, whence $c_i=0$ for all $i$. Thus $A$ is free of rank $q$ over 
$A^G$.

Suppose that $N$ is an $A^G$-module and $M=N\otimes_{A^G}A$ is given a 
$G$-module structure by means of the comodule structure map
$$
\id_N\otimes\,\mu^*:N\otimes_{A^G}A\to N\otimes_{A^G}(A\otimes R)\cong
(N\otimes_{A^G}A)\otimes R.
$$ 
We claim that the assignment $n\mapsto n\otimes1$ yields an isomorphism 
$N\cong M^G$. This amounts to showing that the exactness of the sequence of 
$A^G$-modules 
$$
0\mapright{}A^G\mapright{}A\longmapright5{\mu^*-p_1^*}A\otimes R
$$ 
is preserved under tensoring with $N$ over $A^G$. Now 
$A=A^G\oplus(A^Ga_2+\ldots+A^Ga_q)$ as we assume $a_1=1$. Next, 
$(\mu^*-p_1^*)(A)$ is an $A^G$-submodule of $F$ generated by the elements 
$\mu^*(a_i)-a_i\otimes1$, $\,i=2,\ldots,q$. Note that these elements 
together with $\mu^*(a_1)=1\otimes1$ give a basis for $F'$ over $A$. Then 
the $A$-submodule generated by $\mu^*(a_i)-a_i\otimes1$, $\,i=2,\ldots,q$, 
is a direct summand of $F'$, hence also of $F$. We have seen that $A^G$ is a 
direct summand of $A$, hence also $(\mu^*-p_1^*)(A)$ is a 
direct summand of $F$ as $A^G$-modules. Our claim follows.

Assertion (4) of the theorem is local on $Y$, hence it suffices to consider 
an affine scheme $Y\cong\Spec B$. This is then a special case of what we 
have just proved. Next, taking $N=J$ where $J$ is an ideal of $A^G$, we get 
$JA\cap A^G=J$ since $J\otimes_{A^G}A\cong JA$ by projectivity of $A$ over 
$A^G$. Taking $N=A^G/J$, we deduce that the canonical map $A^G\to(A/JA)^G$ 
is surjective. Suppose that $I$ is a $G$-invariant ideal of $A$. Then 
$\mu^*(I)\subset I\otimes R=IF$. Since on the other hand $\mu^*(A)\subset 
F'$ and $F'$ is a direct summand of $F$, we get $\mu^*(I)\subset F'\cap 
IF=IF'$. Given $a\in I$, we can write therefore the expression $(*)$ with 
$b_i\in I$. As we have seen, this implies $a=\sum b_ia_i$ and $b_i\in A^G$.
Thus $I=I^GA$, which completes the proof of (3). If $\n$ is a maximal ideal 
of $A^G$ then $\n A$ is a maximal $G$-invariant ideal of $A$ by (3). If now 
$\n=\m_x\cap A^G$ then 
$\mu^*(\n A)=(\n\otimes1)\cdot\mu^*(A)\subset\m_x\otimes R$, whence 
$\mu_x^*(\n A)=0$. It follows that $\ker\mu_x^*=\n A$, and $A/\n 
A\cong\mu_x^*(A)$. This proves (2). Assertion (5) follows from (1) since 
$k(X)^G$ is the field of fractions of the ring $A^G$.
\endproof
\proclaim
Corollary 2.3.
If $x$ is a smooth point of $\Xreg$ then $\pi(x)$ is a smooth 
point of $X/G$.
\endproclaim

\Proof.
The local ring $\O_{x,X}$ is a flat extension of $\O_{\pi(x),X/G}$. Since
$\O_{x,X}$ is a regular local ring, so is $\O_{\pi(x),X/G}$ too by
\cite{Ma}, (21.D).
\endproof

\Remark.
Theorem 2.1 can be generalized to the case when $k$ is any field and $X$ is 
a reduced scheme over $k$. In general the stabilizer $G_x$ is a subgroup 
scheme of $G\otimes k(x)$ where $k(x)$ is the residue field of a point 
$x\in X$, and one can define $\Xreg$ to be the set of all points where the 
function $x\mapsto(G\otimes k(x):G_x)$ is locally constant. It can also be 
proved that the morphism $\nu:\Xreg\times G\to\Xreg\times_{\pi(\Xreg)}\Xreg$ 
is finite flat. The assumption that $X$ is reduced is needed to ensure that 
$F'$ is a direct summand of $F$ in the proof of theorem 2.1. Simple 
examples show what can happen without this assumption. Suppose, for 
instance, that $X=\Spec A$ where $A=k[x]$, $\,x^3=0$, and $G$ is the cyclic 
order 2 group with a generator $\sigma$ which acts on $A$ as the 
automorphism sending $x$ to $-x$. Here $A^G=k+kx^2$ (provided $\chr k\ne2$). 
Clearly $A$ is not free over $A^G$. At the same time $X$ contains a single 
point, so that $\Xreg=X$. 
\endremark

\section
3. $G$-linearized modules.

We keep our assumptions on $G$ and $X$ from section 2. Moreover, we assume 
here that $X$ is affine. Let $R=k[G]$, $\,A=k[X]$ and $K=k(X)$. Recall that 
$\M_A$ denotes the category of $(A,G)$-modules. Denote by $\mu^M:M\to M\otimes 
R$ the map that gives $M\in\M_A$ the $R$-comodule structure corresponding to 
the $G$-module structure. In particular, $A$ is an $R$-comodule via the 
map $\mu^*:A\to A\otimes R$ which is the comorphism of $\mu$. We 
may view $M\otimes R$ as a module over $A\otimes R$ in a natural way. The 
compatibility of $A$- and $G$-module structures on $M$ can be 
expressed in terms of the identity $\mu^M(ma)=\mu^M(m)\cdot\mu^*(a)$ where 
$m\in M$ and $a\in A$. We note also that $(S^{-1}M)^G\cong S^{-1}M^G$
for every multiplicatively closed subset $S\subset A^G$.

\proclaim
Lemma 3.1.
Suppose that $X=\Xreg$ and $M\in\M_A$. Then
$$
\eqalign{
M^GA&{}=\{m\in M\mid\mu^M(m)\in(M\otimes1)\cdot\mu^*(A)\}\cr
\noalign{\smallskip}
&{}=\{m\in M\mid m\otimes1\in\mu^M(M)\cdot(A\otimes1)\}.\cr
}
$$
\endproclaim

\Proof.
Let $F=A\otimes R$ and $F'=(A\otimes1)\cdot\mu^*(A)$ as in the proof of 
Theorem 2.1. For every $A$-module $N$ we consider $F_N=N\otimes R$ as an 
$F$-module in a natural way and as an $A$-module by means of the 
homomorphism $p_1^*:A\to F$. Then $F_N\cong N\otimes_AF$. Put 
$F'_N=(N\otimes1)\cdot\mu^*(A)\subset F_N$. Since $F'$ is a 
direct summand of $F$, the canonical map $N\otimes_AF'\to N\otimes_AF$ is a 
split monomorphism of $A$-modules. Clearly its image coincides with $F'_N$. 
Hence $F'_N\cong N\otimes_AF'$.

Suppose that $m\in M$ and $\mu^M(m)\in F'_M$. To prove that $m\in M^GA$ it 
suffices to show that for every $x\in X$ there exists $f\in A^G$ such that 
$f(x)\ne0$ and $mf\in M^GA$. Passing to suitable localizations $A_f$ and 
$M_f$, we may thus assume as in the proof of Theorem 2.1 that $F'$ is a free 
$A$-module with basis elements $\mu^*(a_1),\ldots,\mu^*(a_q)$. Then
$\mu^M(m)=\sum\,(m_i\otimes1)\cdot\mu^*(a_i)$ for some $m_1,\ldots,m_q\in M$.
Let $\epsilon$ and $m^*$ be the counit and the comultiplication in $R$.
Applying $\id_M\otimes\,\epsilon$ to both sides of the equality, we get 
$m=\sum m_ia_i$. Applying $\mu^M\otimes\id_R$ and $\id_M\otimes\,m^*$, we get
$$
\sum\,(\mu^M(m_i)\otimes1)\cdot(\mu^*\otimes\id_R)\mu^*(a_i)= 
\sum\,(m_i\otimes1\otimes1)\cdot(\mu^*\otimes\id_R)\mu^*(a_i)
$$
in $M\otimes R\otimes R$. If $n_1,\ldots,n_q$ are elements of an $A$-module 
$N$ with the property that $\sum\,(n_i\otimes1)\cdot\mu^*(a_i)=0$ in $F_N$ 
then $\sum\,n_i\otimes\mu^*(a_i)=0$ in $N\otimes_AF'$ by the discussion at 
the beginning of the proof, whence $n_i=0$ for all $i$. Now take $N=M\otimes 
R$ with the $A$-module structure given by means of the algebra homomorphism 
$\mu^*:A\to A\otimes R$. Then $A\otimes R$ operates in $N\otimes R\cong 
M\otimes R\otimes R$ by means of the algebra homomorphism 
$\mu^*\otimes\id_R$, and it follows from the displayed equation above that 
$\mu^M(m_i)=m_i\otimes1$, i.e., $m_i\in M^G$ for all $i$. Thus $m\in M^GA$.

Suppose now that $m\in M$ and $m\otimes1=\sum\,\mu^M(m_i)\cdot(b_i\otimes1)$ 
for some elements $m_1,\ldots,m_r\in M$ and $b_1,\ldots,b_r\in A$. If 
$\beta:R\to B$ is the algebra homomorphism corresponding to a point $g\in 
G(B)$ where $B$ is a commutative algebra then, applying 
$\id_M\otimes\,\beta$ to both sides of the equality, we get 
$m\otimes1=\sum\,g(m_i\otimes1)\cdot(b_i\otimes1)$ in $M\otimes B$ (regarded 
as a module over $A\otimes B$). Replacing here $g$ with $g^{-1}$ and 
applying $g$ to both sides of the equality obtained, we get 
$g(m\otimes1)=\sum\,(m_i\otimes1)\cdot g(b_i\otimes1)$. If now $B=R$ and 
$g\in G(R)$ is the point corresponding to the identity homomorphism $R\to 
R$, this can be rewritten as $\mu^M(m)=\sum\,(m_i\otimes1)\cdot\mu^*(b_i)$. 
Thus we have come to the case already considered.  
\endproof

\proclaim
Proposition 3.2.
Suppose that $X=\Xreg$ and denote by $\M'_A$ the full subcategory of $\M_A$ 
consisting of $(A,G)$-modules $M$ such that $M=M^GA$. Then:

\item
(1) $\M'_A$ is closed under taking submodules and factor modules.

\item
(2) The functor $M\mapsto M^G$ is an equivalence between $\M'_A$ and the 
category of $A^G$-modules. The inverse functor is $N\mapsto 
N\otimes_{A^G}A$.

\item
(3) If $M\in\M'_A$ is projective of rank $r$ as an $A$-module then $M^G$ is 
projective of rank $r$ as an $A^G$-module. 

\endproclaim

\Proof. (1) We use the same notations as in the preceding lemma. Clearly, 
$F_N/F'_N\cong N\otimes_AF/F'$ for every $A$-module $N$. Let $N$ be an 
$(A,G)$-submodule of $M\in\M'_A$. Since $F/F'$ is a projective $A$-module, 
the canonical map $N\otimes_AF/F'\to M\otimes_AF/F'$ is injective. It 
follows from the commutative diagram 
$$
\diagram{
0&\mapright{}&F'_N&\mapright{}&F_N&\mapright{}&N\otimes_AF/F'&\mapright{}&0\cr
\diagramskip
&&\downarrow&&\downarrow&&\downarrow&&\cr
\diagramskip
0&\mapright{}&F'_M&\mapright{}&F_M&\mapright{}&M\otimes_AF/F'&\mapright{}&0\cr
}
$$
that $F'_N=F_N\cap F'_M$. Now $\mu^M(N)\subset F_N$ and $\mu^M(M)\subset F'_M$,
whence $\mu^M(N)\subset F'_N$. Hence $N\in\M'_A$ by Lemma 3.1.
The assertion about factor modules is obvious.

(2) If $N$ is an $A^G$-module and $M=N\otimes_{A^G}A$ then $N\cong M^G$ as 
we have seen in the proof of Theorem 2.1.
Conversely, suppose that $M\in\M'_A$ and $M_1=M^G\otimes_{A^G}A$. The 
canonical map $\phi:M_1\to M$ is a morphism in $\M_A$. It is surjective by 
the definition of $\M'_A$. Now $\ker\phi\in\M'_A$ by (1), hence $\ker\phi$ 
is generated over $A$ by $G$-invariant elements. But we have proved already 
that $M_1^G=M^G\otimes1$. Since $\phi$ is injective on $M_1^G$, we get 
$\ker\phi=0$, i.e., $\phi$ is an isomorphism.

(3) follows from \cite{Bou}, Ch.~I, \S3, Prop.~12 and Ch.~II, \S5, Prop.~4.
\endproof

Suppose that $M\in\M_A$ is finitely generated over $A$. Put
$$
\rk M=\dim_KM\otimes_AK\qquad{\fam0and}\qquad 
M(x)=M/M\m_x\quad{\fam0for\ }x\in X  
$$
where $\m_x$ is the maximal ideal of $A$ corresponding to $x$.
Then $\dim M(x)\ge\rk M$ for all $x\in X$. By \cite{Bou}, Ch.~2, \S3, 
Prop.~7 and \S5, Corollary to Prop.~2, the set 
$$
X_M=\{x\in X\mid\dim M(x)=\rk M\}
$$ 
is open in $X$ and consists precisely of those $x$ for which 
$M_{\m_x}$ is a free $A_{\m_x}$-module. If $\m_x$ is stable under $G$ then
$\mu^M(M\m_x)=\mu^M(M)\cdot\mu^*(\m_x)\subset M\m_x\otimes R$ since 
$\mu^*(\m_x)\subset\m_x\otimes R$. In general, applying this observation to 
the action of $G_x$, we see that $M\m_x$ is stable under $G_x$, and so $G_x$ 
operates in $M(x)$. Put 
$$ 
s(M)=\min_{x\in\Xreg\cap X_M}\dim M(x)^{G_x},
$$\removelastskip
$$ 
\Mreg=\{x\in\Xreg\mid\dim M(x)=\rk M\quad{\fam0and}\quad
\dim M(x)^{G_x}=s(M)\}.  
$$
We call $\Mreg$ the set of $M$-regular points in $X$.

\proclaim
Theorem 3.3. 
{\rm(1)} $\Mreg$ is a $G$-invariant open subset of $X$. 

\item
(2) For all $x\in\Mreg$ the canonical map $M\to M(x)$ induces a surjection 
$M^G\to M(x)^{G_x}$.

\item
(3) If $\Mreg=X$ then the map $M^G\otimes_{A^G}A\to M$ given by $\,m\otimes 
a\mapsto ma$ is a split monomorphism of $A$-modules and $M^G$ is 
projective of rank $s(M)$ over $A^G$.

\item
(4) $\dim_{K^G}(M\otimes_AK)^G=s(M)$.

\endproclaim

\Proof.
If $x\in X$ and $g\in G(k)$ then $g$ induces a linear isomorphism $M(x)\to 
M(xg)$ compatible with the actions of stabilizers. Hence $M(x)^{G_x}\cong 
M(xg)^{G_{xg}}$. It follows that $X_M$ and $\Mreg$ are stable under the 
action of $G(k)$. Then $X_M$ is a $G$-invariant open subset, and so is 
$\Xreg\cap X_M$ too. Localizing if necessary, we may assume that $X=\Xreg$ 
and $M_{\m_x}$ is free over $A_{\m_x}$ for all $x$. Then $M$ is a projective 
$A$-module. 

We are going to apply Lemma 1.3 in which we take $F=M\otimes R$ with the 
$A$-module structure obtained again via $p_1^*:A\to A\otimes R$. Take 
$F''=M\otimes1$, which is clearly a direct summand of $F$. Put 
$F'=\mu^M(M)\cdot(A\otimes1)$, and let $I'(x)$, $I''(x)$ be the images, 
respectively, of $F'/F'\m_x$ and $F''/F''\m_x$ in
$F/F\m_x\cong M(x)\otimes R$. We have $I''(x)\cong M(x)\otimes1$. If 
$\mu^M_x$ denotes the composite 
$$ M\longmapright2{\mu^M}M\otimes R
\longmapright6{\can\otimes\,\id_R}M(x)\otimes R, 
$$ 
then $I'(x)=\mu^M_x(M)$. Consider two $G$-module structures on $M\otimes R$: 
the first one is the tensor product of the given $G$-module structure on $M$ 
and the left regular $G$-module structure on $R$; the second one is the 
tensor product of the trivial $G$-module structure on $M$ and the right 
regular $G$-module structure on $R$. The map $\mu^M$ is $G$-equivariant with 
respect to the second structure and is a bijection of $M$ onto the subspace 
$(M\otimes R)^G$ of $G$-invariant elements with respect to the first 
structure (see \cite{Ja}, Part I, 3.7, (5) and (6); however, we interchanged 
the left and right regular $G$-module structures). These two structures on 
$M\otimes R$ induce a $G_x$-module and a $G$-module structures on 
$M(x)\otimes R$. We get $$ \mu^M_x(M)\subset(M(x)\otimes 
R)^{G_x}=\ind_{G_x}^GM(x).  $$ Since $\mu^M_x$ is $G$-equivariant, 
$\mu^M_x(M)$ is a $G$-submodule of the induced module. Furthermore, 
$\mu^M_x(ma)=\mu^M_x(m)\cdot\mu_x^*(a)$ for $m\in M$, $a\in A$, whence 
$\mu^M_x(M)$ is stable under the action of $\mu_x^*(A)\subset R$. By 
Proposition 1.1 $\mu_x^*(A)=k[G_x\backslash G]$.  The canonical 
$G_x$-equivariant map $\phi:\ind_{G_x}^GM(x)\to M(x)$ is the restriction of 
the map $\id\otimes\,\epsilon:M(x)\otimes R\to M(x)$.  Hence the composite 
$\phi\circ\mu^M_x$ coincides with the canonical projection $M\to M(x)$, it 
is therefore surjective. It follows that the inclusion 
$\iota:\mu^M_x(M)\hookrightarrow\ind_{G_x}^GM(x)$ corresponds under the 
equivalence of the Imprimitivity Theorem to a surjective map of 
$G_x$-modules. Then $\iota$ must itself be surjective, i.e., 
$\mu^M_x(M)=\ind_{G_x}^GM(x)$. In particular, 
$$
\dim I'(x)=\dim\ind_{G_x}^GM(x)=(G:G_x)\cdot\dim M(x)=q(X)\rk(M), 
$$
which does not depend on $x$. By Lemma 1.2 $F'$ is a 
direct summand of $F$. Thus the hypotheses of Lemma 1.3 are fulfilled. By our 
previous description $I'(x)\cap I''(x)=(M(x)\otimes1)^{G_x}\cong 
M(x)^{G_x}$. Hence $s=s(M)$ and $\Mreg$ is the set of rational points of the 
open subset $U\subset\Spec A$ defined in Lemma 1.3. Thus $\Mreg$ is open in 
$X$. It is $G$-invariant by observation at the beginning of the proof.  
Therefore we may localize further and assume $X=\Mreg$. By Lemma 3.1 $F'\cap 
F''=M^GA\otimes1$. Then (2) of Lemma 1.3 implies (2) of the theorem. Lemma 
1.3 ensures that the $A$-module $M^GA$ is projective of rank $s(M)$ and is a 
direct summand of the $A$-module $M$. By Proposition 3.2 $M^G$ is a 
projective $A^G$-module of rank $s(M)$ and $M^G\otimes_{A^G}A$ is mapped 
isomorphically onto $M^GA$. This proves (3). Assertion (4) is immediate 
since $(M\otimes_AK)^G\cong M^G\otimes_{A^G}K^G$.  \endproof

Let $M,N\in\M_A$ and $P=\Hom_A(N,M)$. For every finite dimensional 
commutative algebra $B$ we have
$$
P\otimes B\cong\Hom_{A\otimes B}(N\otimes B,\,M\otimes B).
$$                                
If $g\in G(B)$ and $\xi\in P\otimes B$ then we put 
$g_P(\xi)=g_M\circ\xi\circ g_N^{-1}$ where $g_M$ and $g_N$ are the operators 
on $M\otimes B$ and $N\otimes B$, respectively, corresponding to $g$. In 
this way we obtain a group action of $G(B)$ on $P\otimes B$ which is natural 
in $B$. If $B$ is infinite dimensional then each point $g\in G(B)$ still 
belongs to $G(B')$ where $B'\subset B$ is a finite dimensional subalgebra. 
Indeed, we can take $B'$ to be the image of the algebra homomorphism $R\to 
B$ corresponding to $g$. Extend the action of $g$ in $P\otimes B'$ by 
$B$-linearity to the action in $P\otimes B$. If $g,h\in G(B)$ are two points 
then there exists a finite dimensional subalgebra $B'$ such that $G(B')$ 
contains both of them. It follows that $(gh)_P=g_Ph_P$. Thus $P$ is equipped 
with a $G$-module structure, which is clearly compatible with the $A$-module 
structure, i.e., $P\in\M_A$.

Let $V$ be a $G$-module. Then $V\otimes A$, considered with the natural 
$A$-module structure and the tensor product $G$-module structure, is an 
object of $\M_A$. Hence so is $\Hom(V,M)\cong\Hom_A(V\otimes A,M)$ too.

Suppose that $M,N\in\M_A$ are finitely generated over $A$ and $V$ a finite 
dimensional $G$-module. Put
$$
\trialign{ 
s(N,M)={}&\min_{x\in\Xreg\cap X_M\cap X_N}{}
&\dim\Hom_{G_x}\bigl(N(x),\,M(x)\bigr),\cr
\noalign{\smallskip}
s(V,M)={}&\min_{x\in\Xreg\cap X_M}{}&\dim\Hom_{G_x}\bigl(V,\,M(x)\bigr).\cr
}
$$

\proclaim 
Corollary 3.4. 
{\rm(1)} $\dim_{K^G}\Hom_{(K,G)}(N\otimes_AK,\,M\otimes_AK)=s(N,M)$.

\item
(2) $\dim_{K^G}\Hom_G(V,\,M\otimes_AK)=s(V,M)$.

\item
(3) $\dim_{K^G}\soc_G(M\otimes_AK)=\sum\,s(V,M)\dim V$, the sum over 
isomorphism classes of irreducible $G$-modules $V$.

\endproclaim

\Proof. (1) Let $P=\Hom_A(N,M)$. Then
$P\otimes_AK\cong\Hom_K(N\otimes_AK,\,M\otimes_AK)$ and
$P(x)\cong\Hom\bigl(N(x),M(x)\bigr)$ for all $x\in X_N$. Since $X_M\cap 
X_N\subset X_P$, we get $s(P)=s(N,M)$. Apply Theorem 3.3(4) to the 
$(A,G)$-module $P$, noting that 
$$\Hom_{(K,G)}(N\otimes_AK,\,M\otimes_AK)\cong(P\otimes_AK)^G.$$

Assertion (2) is a special case of (1) with $N=V\otimes A$.
Next, the $G$-socle $\soc_G(M\otimes_AK)$ is a direct sum of isotypic 
components $I_V$ corresponding to irreducible $G$-modules $V$. Since 
$\End_GV\cong k$, we have $I_V\cong\Hom_G(V,\,M\otimes_AK)\otimes V$, and 
(3) follows from (2).  
\endproof

We continue to assume that $M\in\M_A$ is finitely generated over $A$. 

\proclaim 
Lemma 3.5.
If $x\in\Xreg\cap X_M$ and $\n=\m_x\cap A^G$ then 
$M/M\n\cong\ind_{G_x}^GM(x)$ as $G$-modules. The restriction of the 
canonical map $M/M\n\to M(x)$ yields a linear isomorphism 
$(M/M\n)^G\cong M(x)^{G_x}$.
\endproclaim 

\Proof.
Since $\Xreg\cap X_M$ is a $G$-invariant open subset, there exists $f\in 
A^G\setminus\n$ such that $A_f$ is free of rank $q(X)$ over $A_f^G$ and 
$M_f$ is free over $A_f$. Then $M_f$ is free of rank $q(X)\rk(M)$ over 
$A_f^G$, and so $M/M\n\cong M_f/M_f\n$ has dimension $q(X)\rk(M)$. In the 
proof of Theorem 3.3 we constructed a surjective $G$-module homomorphism
$\mu^M_x:M\to\ind_{G_x}^GM(x)$. Since 
$\mu^M(M\n)=\mu^M(M)\cdot(\n\otimes1)\subset M\m_x\otimes R$, we see that 
$M\n\subset\ker\mu^M_x$. Comparing dimensions, we conclude that 
$M\n=\ker\mu^M_x$. The final assertion is a special case of the Frobenius 
reciprocity. 
\endproof

Put $\Minj=\{x\in\Xreg\cap X_M\mid M(x)$ is an injective $G_x$-module$\}$.

\proclaim 
Theorem 3.6. {\rm(1)} $\Minj$ is open in $X$ and consists precisely of those 
points $x\in\Xreg\cap X_M$ for which there exists $f\in A^G$ such that 
$f(x)\ne0$ and $M_f$ is a projective $(A_f,G)$-module.

\item
(2) $\Minj$ is nonempty if and only if $M\otimes_AK$ is a projective 
$(K,G)$-module.

\item
(3) The induced $G$-modules $\ind_{G_x}^GM(x)$ corresponding to points 
$x\in\Minj$ are all isomorphic to each other.

\item
(4) For every $x\in\Minj$ there exist $f\in A^G$ and a $G$-submodule 
$V\subset M$ such that $f(x)\ne0$, $\,V\cong\ind_{G_x}^GM(x)$ and the 
linear map $V\otimes A_f^G\to M_f$ given by the rule $m\otimes a\mapsto ma$ 
is bijective. In particular, $M_f$ is an injective $G$-module.

\item
(5) If $\Minj$ is nonempty then $M\otimes_AK\cong V\otimes K^G$ as 
$(K^G,G)$-modules. In particular, $M\otimes_AK$ is an injective $G$-module.

\endproclaim 

\Proof.
Let $x\in\Minj$. The $G$-module $V_x=\ind_{G_x}^GM(x)$ is injective by 
\cite{Ja}, Part~I, 3.9, hence it is also projective. Then, in view of Lemma 
3.5, there exists a $G$-submodule $V\subset M$ such that $M=V\oplus M\n$ 
where $\n=\m_x\cap A^G$. Clearly $V\cong\nobreak V_x$. Define a homomorphism
of $A^G$-modules $\phi:\nobreak V\otimes\nobreak A^G\to M$ by the formula
$\phi(m\otimes a)=ma$ for $m\in V$, $a\in A^G$. Since
$x\in\nobreak\Xreg\cap\nobreak X_M$, the algebra $A_\n$ is free over 
$A_\n^G$ and $M_\n$ is free over $A_\n$. Then $M_\n$ is free over $A_\n^G$, 
and it follows that $\phi$ induces an isomorphism $V\otimes\nobreak 
A_\n^G\to M_\n$. There exists $f\in A^G\setminus\n$ such that 
$\phi_f:V\otimes A_f^G\to M_f$ is an isomorphism (see \cite{Bou}, Ch.~3, 
\S5, Prop.~2).  As a $G$-module, $M_f$ is a direct sum of a family of copies 
of the $G$-module $V$, and so it is injective. This proves (4), and (5) is 
an immediate consequence. 

If $\n'$ is a maximal ideal of $A^G$ such that 
$f\notin\n'$ then $A_f^G=k\oplus A_f^G\n'$, whence $M_f=V\oplus M_f\n'$, and 
$M/M\n'\cong M_f/M_f\n'\cong V$ as $G$-modules. Again by Lemma 3.5 
$\ind_{G_y}^GM(y)\cong V$ for all $y\in\Xreg\cap X_M$ such that $f(y)\ne0$. If 
$x'\in\Minj$ is another point then, similarly, $\ind_{G_y}^GM(y)\cong 
V_{x'}$ for all $y$ in a nonempty open subset of $X$. Since $X$ is 
irreducible, we conclude $V_x\cong V_{x'}$, whence (3).

Suppose that $N\in\M_A$ is free of finite rank over $A$ and $\phi:N\to M$ is 
an epimorphisms in $\M_A$. Let $\phi_x:N(x)\to M(x)$ be the epimorphism of 
$G_x$-modules obtained from $\phi$ by reduction modulo $\m_x$. Since 
$M(x)$ is projective, there exists a $G_x$-module homomorphism 
$\psi_x:M(x)\to N(x)$ such that $\phi_x\circ\psi_x=\id_{M(x)}$. Put 
$P=\Hom_A(M,N)$. Then $X_M\subset X_P$ and for $x\in X_M$ there is an 
isomorphism of $G_x$-modules
$P(x)\cong\Hom\bigl(M(x),N(x)\bigr)\cong M(x)^*\otimes N(x)$.
Since $\H(G_x)$ is a Frobenius algebra, the $G_x$-module $M(x)^*$ is 
injective. By \cite{Ja}, Part I, 3.10 $P(x)$ is also injective. As we know 
already, there exists a $G$-submodule $W\subset P$ such that $P=W\oplus 
P\n$. By Lemma 3.5 the restriction of the canonical map $W\iso P/P\n\to P(x)$ 
yields a linear isomorphism $W^G\iso P(x)^{G_x}$. Let $\psi\in 
W^G\subset P^G$ be the element corresponding to $\psi_x\in P(x)^{G_x}$. We 
may regard $\psi$ as a morphism $M\to N$ in $\M_A$ whose reduction modulo 
$\m_x$ is $\psi_x$. Then $\gamma=\phi\circ\psi$ is an $\M_A$-endomorphism 
of $M$ whose reduction modulo $\m_x$ is the identity transformation of 
$M(x)$.  Let $U$ be the set of those $y\in X_M$ for which the reduction of 
$\gamma$ modulo $\m_y$ is invertible. By \cite{Bou}, Ch.~2, \S3, Corollary 
to Prop.~6 and \S5, Prop.~2 $U$ is open and consists precisely of those 
$y\in X_M$ for which $\gamma_{\m_y}:M_{\m_y}\to M_{\m_y}$ is bijective. 
Since $\gamma$ is $G(k)$-equivariant, $U$ is $G$-invariant. As 
$U=\pi^{-1}\bigl(\pi(U)\bigr)$ and $\pi(U)$ is an open neighbourhood of 
$\pi(x)$ in $X/G$, there exists $f\in A^G$ such that $x\in X_f\subset U$ 
where $X_f=\{y\in X\mid f(y)\ne0\}$. Then $\gamma_f:M_f\to M_f$ is 
bijective, and therefore $N_f=\ker\phi_f\oplus\im\psi_f$. In other words, 
$\phi_f:N_f\to M_f$ is a split epimorphism in $\M_{A_f}$. Since $M$ is 
finitely generated over $A$, it is an epimorphic image of a finitely 
generated free $\smashprod{A}$-module. We can take the latter to be our $N$.  
We see that $M_f$ is a direct summand of a free $\smashprod{A_f}$-module for 
a suitable $f$.

Conversely, suppose that $x$ is any point in $\Xreg\cap X_M\cap X_f$ such 
that $M_f$ is a projective $\smashprod{A_f}$-module. If $\n=\m_x\cap A^G$ 
then $M/M\n\cong M_f/M_f\n$ is a projective $\smashprod{(A/A\n)}$-module. By 
Theorem 2.1(2) $A/A\n\cong k[G_x\backslash G]$. 
The $G_x$-module $M(x)$ corresponds to $M/M\n$ under the category 
equivalence of the Imprimitivity Theorem. It is therefore projective, hence 
injective. We get (1).

Suppose that $M\otimes_AK$ is projective in $\M_K$. We want to show that 
$M_f$ is projective in $\M_{A_f}$ for a suitable $0\ne f\in A^G$ and then 
apply (1). Since $X_M$ is open and $G$-invariant, we may assume that $X_M=X$ 
passing at the very beginning to a suitable localization of $A$. Then $M$ is 
projective as an $A$-module. Let $\phi:N\to M$ be an epimorphism in $\M_A$ 
with $N$ a free $\smashprod{A}$-module. It extends to an epimorphism 
$\phi_K:N\otimes_AK\to M\otimes_AK$ in $\M_K$. By our assumptions the latter 
admits a splitting $\psi:M\otimes_AK\to N\otimes_AK$ in $\M_K$. Since $N$ 
is free over $A$, the localizations $N_f$ are identified with their images 
in $N\otimes_AK$, and the same is valid for $M$. Since $M$ is finitely 
generated over $A$, hence also over $A^G$, we have $\psi(M)\subset N_f$ for 
a suitable $f$. Then $\psi(M_f)=N_f$, which means that $\phi_f:N_f\to M_f$ 
is a split epimorphism in $\M_{A_f}$. That completes the proof of (2).  
\endproof

\section
4. Actions with linearly reductive stabilizers.

We weaken our assumptions for this section considerably. In the next 
proposition $G$ is any affine group scheme over $k$, not necessarily finite, 
and $A$ is any $G$-algebra, not necessarily commutative. What we prove is a 
special case of results due to Doi \cite{Doi} obtained in the context of 
coactions of Hopf algebras.

\proclaim
Proposition 4.1. The following properties of a $G$-algebra $A$ are 
equivalent:

\item
(1) All objects $M\in\M_A$ are injective $G$-modules.

\item
(2) There exist an injective $G$-module $Q$ and a homomorphism
of $G$-modules $\psi:Q\to A$ such that $1\in\psi(Q^G)$.

\item
(3) There exists a $G$-module homomorphism $\phi:k[G]\to A$ {\rm(}where 
$k[G]$ is given the left regular $G$-module structure{\rm)} such that 
$\phi(1)=1$.

\item
(4) There are linear maps $\Phi_M:M\to M^G$, defined for each $M\in\M_A$, 
which are natural in $M$ and satisfy $\Phi_M(m)=m$ for all $m\in M^G$.

\item
(5) The functor $M\mapsto M^G$ is exact on $\M_A$.

\item
(6) $A$ is a projective $(A,G)$-module.

If $A$ is commutative they are equivalent to another property:

\item
(7) Every object $M\in\M_A$ which is finitely generated and projective as an
$A$-module is projective in $\M_A$.

\endproclaim

\Proof.
$(1)\Rightarrow(2)$.
By (1) $A$ is an injective $G$-module. So we can take $Q=A$ and $\psi=\id_A$.

$(2)\Rightarrow(3)$.
Let $q\in Q^G$ be an element such that $\psi(q)=1$. By injectivity of $Q$ 
the $G$-module homomorphism $k\to Q$ sending $1$ to $q$ extends to a 
homomorphism $k[G]\to Q$. Composing the latter with $\psi$, we get $\phi$.

$(3)\Rightarrow(4)$.
Define $\Phi_M$ as the composite map
$$
M\longmapright2{\mu^M}M\otimes k[G]
\longmapright4{\id_M\otimes\,\phi}M\otimes A\longmapright2{}M
$$
where the last map is afforded by the $A$-module structure on $M$. Recall 
that $\mu^M(M)=(M\otimes k[G])^G$. Since the two final maps in the 
decomposition of $\Phi_M$ are $G$-equivariant, we get $\Phi_M(M)\subset 
M^G$. If $m\in M^G$ then $\mu^M(m)=m\otimes1$, whence $\Phi_M(m)=m$. That the 
maps $\Phi_M$ are natural in $M$ is clear.

$(4)\Rightarrow(5)$.
The fixed point functor is clearly left exact. Suppose that $\xi:M\to N$ is 
an epimorphism in $\M_A$. Given $n\in N^G$, take $m\in M$ such that 
$\xi(m)=n$. Then $\Phi_M(m)\in M^G$ and 
$\xi\bigl(\Phi_M(m)\bigr)=\Phi_N\bigl(\xi(m)\bigr)=n$. Thus $\xi$ induces a 
surjection $M^G\to N^G$.

$(5)\Rightarrow(1)$.
We may view $\Hom(V,M)\cong V^*\otimes M$ for each finite dimensional 
$G$-module $V$ as an $(A,G)$-module taking the tensor product of $G$-module 
structures and letting $A$ operate on the second tensorand. If $W\subset V$ 
is a $G$-submodule then we have an epimorphism $\Hom(V,M)\to\Hom(W,M)$ in 
$\M_A$. Applying the fixed point functor, we deduce the surjectivity of the 
canonical map $\Hom_G(V,M)\to\Hom_G(W,M)$. Since all $G$-modules are locally 
finite dimensional, this gives the injectivity of $M$.

$(5)\Leftrightarrow(6)$.
Every morphism $A\to M$ in $\M_A$ is given by the rule $a\mapsto ma$ where 
$m\in M^G$. Hence $\Hom_{(A,G)}(A,M)\cong M^G$. Note that the projectivity 
of $A$ in $\M_A$ means that the functor $M\mapsto\Hom_{(A,G)}(A,M)$ is 
exact.

$(5)\Rightarrow(7)$.
If $M$, $N$ are $A$-modules, $P=\Hom_A(M,N)$ and $B$ a commutative algebra 
then the canonical map $P\otimes B\to\Hom_{A\otimes B}(M\otimes B,\,N\otimes 
B)$ is bijective when $M$ is free of finite rank, hence also when $M$ is 
finitely presented. If $M,N\in\M_A$ and $M$ is finitely presented as an 
$A$-module then $G(B)$ operates in $P\otimes B$, naturally in $B$. 
This gives $P$ a $G$-module structure. Assuming $A$ to be commutative, we 
have $P\in\nobreak\M_A$. If, moreover, $M$ is projective as an $A$-module 
then every epimorphism $N\to N'$ in $\M_A$ induces an epimorphism 
$\Hom_A(M,N)\to\Hom_A(M,N')$. Applying the fixed point functor, we deduce 
the surjectivity of the map 
$\Hom_{(A,G)}(M,N)\to\Hom_{(A,G)}(M,N')$.

$(7)\Rightarrow(6)$ is obvious.
\endproof

Let $X$ be an arbitrary scheme over $k$, and $G$ a finite group scheme 
operating on $X$ from the right. We still need the assumption that $X$ can 
be covered by $G$-invariant affine open subschemes. We say that the 
stabilizer $G_x$ of a point $x\in X(k)$ is {\it linearly reductive} if all 
$G_x$-modules are completely reducible. This is equivalent to the 
semisimplicity of the Hopf algebra $\H(G_x)$. By \cite{DG}, Ch.~IV, \S3, 3.6 
$G_x$ is linearly reductive if and only if its identity component $G_x^0$ 
is diagonalizable and the index $(G_x:G_x^0)$ is prime to $p=\chr k$ when 
$p>0$. Put 
$$ 
\Xlr=\{x\in X(k)\mid G_x{\fam0\ is\ linearly\ reductive}\}.
$$

\proclaim
Theorem 4.2.
The set $\Xlr$ consists precisely of those $x\in X(k)$ which are contained 
in a $G$-invariant affine open subscheme $U\subset X$ such that $k[U]$ is an 
injective $G$-module. In particular, $\Xlr$ is the set of rational points of 
an open $G$-invariant subscheme of $X$. If, moreover, $X$ is an algebraic 
variety, then:

\item 
(1) The condition $\Xlr\ne\nothing$ is equivalent to each of the two below:

\subitem 
(a) $k(X)$ is an injective $G$-module.

\subitem 
(b) The smash product algebra $\smashprod{k(X)}$ is semisimple. 

\item 
(2) For every $x\in\Xreg\cap\Xlr$ there exist a $G$-invariant affine open 
neighbourhood $U$ of $x$ and a $G$-submodule $V\subset k[U]$ such that 
$V\cong\ind_{G_x}^Gk$ and the map $V\otimes k[U]^G\to k[U]$ given by the 
multiplication in $k[U]$ is bijective.

\item 
(3) If $\Xlr$ is nonempty then there exists a $G$-submodule $V\subset k[X]$ 
such that $V\cong\ind_{G_x}^Gk$ for all $x\in\Xreg\cap\Xlr$ and the map 
$V\otimes k(X)^G\to k(X)$ given by the multiplication in $k(X)$ is 
bijective.

\endproclaim

\Proof.
We may assume that $X$ is affine. Put $A=k[X]$. Given $x\in X(k)$, the orbit 
morphism $\mu_x:G\to X$ determines a $G$-equivariant homomorphism of algebras
$\mu_x^*:A\to k[G]$ whose image is $k[G_x\backslash G]$ by Proposition 1.1.
If $x\in\Xlr$ then all $G_x$-modules are injective. Hence 
$k[G_x\backslash G]=\ind_{G_x}^Gk$ is an injective $G$-module by \cite{Ja}, 
Part~I, 3.9. As it is also projective, there is a $G$-submodule $V\subset A$ 
mapped isomorphically onto $k[G_x\backslash G]$ under $\mu_x^*$. Take 
$f\in V$ such that $\mu_x^*(f)=1$. Then $f\in A^G$, and the map $V\to A_f$, 
$\,v\mapsto vf^{-1}$, is a $G$-module homomorphism under which $f\mapsto1$.
Since $V$ is an injective $G$-module, so is $A_f$ too by implication 
$(2)\Rightarrow(1)$ of Proposition 4.1. Furthermore, $f(x)=1$ since 
$\ker\mu_x^*\subset\m_x$. Thus $\Spec A_f$ is the required open 
neighbourhood of $x$. 

Conversely, suppose that $x\in U(k)$ where $U\subset X$ is a $G$-invariant 
affine open subscheme such that $k[U]$ is an injective $G$-module. As 
$\mu_x$ factors through $U$, it induces a $G$-equivariant algebra 
homomorphism $k[U]\to k[G_x\backslash G]$. Hence we may view 
$k[G_x\backslash G]$ as a $(k[U],G)$-module. By Proposition 4.1
$k[G_x\backslash G]$ is an injective $G$-module. It is then a direct summand 
of $k[G]$. By \cite{Ja}, Part I, 4.12 $k[G_x\backslash G]$ is also an 
injective $G_x$-module. Now $k$ is a direct summand of $k[G_x\backslash G]$ 
as a $G_x$-module. Hence $k$ is an injective $G_x$-module. Then all 
$G_x$-modules are injective, which implies that all $G_x$-modules are 
completely reducible.

Suppose that $X$ is an algebraic variety and $K=k(X)$. Apply Theorem 3.6 
to the $(A,G)$-module $A$. Noting that $\Xlr$ is precisely the set of points 
$x$ for which $A/\m_x\cong k$ is an injective $G_x$-module, we get 
assertions (2) and (3) of Theorem 4.2. Furthermore, $\Xlr$ is nonempty if 
and only if $K$ is a projective $(K,G)$-module. By Proposition 4.1 this 
is equivalent to $K$ being an injective $G$-module. This is equivalent 
also to the condition that every $M\in\M_K$ of finite dimension over 
$K$ is projective in $\M_K$. This means, in particular, that 
all ideals of the algebra $\smashprod{K}$ are projective, which is 
equivalent to condition (b).
\endproof

\section
5. Invariants of restricted Lie algebras.

Suppose that $\chr k=p>0$. Let $X$ be an affine algebraic variety, and $\g$ 
a $p$-Lie algebra over $k$. Put $A=k[X]$ and $K=k(X)$. Define an {\it 
action} of $\g$ on $X$ to be a homomorphism of $p$-Lie algebras 
$\rho:\g\to\Der A$ into the derivation algebra of $A$. Define 
$\g_x\subset\g$ to be the stabilizer of the maximal ideal $\m_x$ of $A$ 
corresponding to a point $x\in X$. Since $\rho(\g)(\m_x^2)\subset\m_x$, we 
have a linear map $\g\to T_xX=(\m_x/\m_x^2)^*$ whose kernel is precisely 
$\g_x$. Hence $\codim_\g\g_x\le\dim T_xX$. Since the dimensions of tangent 
spaces are bounded, it is meaningful to define
$$
c_\g(X)=\max_{x\in X}\codim_\g\g_x,
$$\removelastskip 
$$
\greg=\{x\in X\mid\codim_\g\g_x=c_\g(X)\}.
$$
If $\dim\g<\infty$, there is a finite group scheme of height one $G=\G(\g)$ 
associated with $\g$ (see \cite{DG}, Ch.~II, \S7, 3.9). One has $k[G]\cong 
u(\g)^*$ where $u(\g)$ is the restricted universal enveloping algebra of 
$\g$. The action of $\g$ on $X$ corresponds to a group action of $G$ 
according to \cite{DG}, Ch.~II, \S7, 3.10. Furthermore, $G_x\cong\G(\g_x)$, 
so that $k[G_x\backslash G]\cong\Hom_{u(\g_x)}\bigl(u(\g),k\bigr)$ and 
$(G:G_x)=p^{\codim_\g\g_x}$ for all $x\in X$. It follows then that 
$\greg=\Xreg$. However, we want to extend Theorem 2.1 to the case of 
infinite dimensional $\g$.

The Lie algebra $\Der A$ has a natural $A$-module structure. Given $f\in A$ 
and $D,D'\in\Der A$, we have 
$$ 
[fD,D']=f\,[D,D']-D'(f)D, 
$$\removelastskip 
$$ 
(fD)^p=f^pD^p+(fD)^{p-1}(f)D.
$$
The first formula is easily checked straightforwardly. The second one is 
proved by Hochschild \cite{Hoch2}, Lemma 1. It follows that the 
$A$-submodule $L=A\cdot\rho(\g)$ is also a $p$-Lie subalgebra of $\Der A$. 
Define a linear map $$d:A\to L_A^*=\Hom_A(L,A)$$ by the rule $(df)(D)=D(f)$ 
for $f\in A$, $D\in L$.

\proclaim
Lemma 5.1. {\rm(1)} $c_\g(X)=c_L(X)$ and $\greg=\Lreg$.

\item
(2) $\Lreg$ is open in $X$ and consists precisely of those $x\in X$ for 
which $L_{\m_x}$ is a free $A_{\m_x}$-module and $L_A^*=dA+\m_xL_A^*$.

\item
(3) If $\Lreg=X$ then $L$ is projective of rank $c_L(X)$ over $A$
and $L_A^*=A\cdot dA$. 

\endproclaim

\Proof. If $L_x$ is the stabilizer of $\m_x$ in $L$ then $\m_xL\subset L_x$, 
whence $L=\rho(\g)+L_x$. It follows that $L/L_x\cong\g/\g_x$, and (1) is 
immediate.

Take a finite system of generators $a_1,\ldots,a_n$ of the algebra $A$ and 
define a homomorphism of $A$-modules $\phi:L\to F$, where $F=A^n$, by the 
rule $\phi(D)=(Da_1,\ldots,Da_n)$ for $D\in L$. Since each 
derivation of $A$ is determined by its values on generators, we have 
$\ker\phi=0$. Hence $L\cong F'=\phi(L)\subset F$. Denote by $I(x)$ the image 
of the map $L/\m_xL\cong F'/\m_xF'\to F/\m_xF$ induced by $\phi$. If $D\in 
L$ then $D\in L_x$ if and only if $D(A)\subset\m_x$ (as $A=k+\m_x$), if and 
only if $Da_i\in\m_x$ for all $i$.  It follows that $I(x)\cong L/L_x$.  
Applying Lemma 1.2, we see that $\Lreg$ is open and coincides with the set 
of points $x\in X$ for which $F'_{\m_x}$ is a direct summand of the 
$A_{\m_x}$-module $F_{\m_x}$, i.e., $\phi_{\m_x}:L_{\m_x}\to F_{\m_x}$ is a 
split monomorphism of $A_{\m_x}$-modules. Assuming that the 
$A_{\m_x}$-module $L_{\m_x}$ is free, $\phi_{\m_x}$ splits if and only if 
the localization at $\m_x$ of the dual homomorphism 
$\phi^*:\Hom_A(F,A)\to\Hom_A(L,A)$ is surjective. By Nakayama's lemma this 
is equivalent to the equality $L_A^*=N+\m_xL_A^*$ where $N$ is the image of 
$\phi^*$. Clearly, $N$ is the $A$-submodule of $L_A^*$ generated by 
$da_1,\ldots,da_n$. Since $d(ab)=a\cdot db+b\cdot da$ for $a,b\in A$, we 
have $N=A\cdot dA$. Then $N+\m_xL_A^*=dA+\m_xL_A^*$ since $A=k+\m_x$. We get 
(2).

If $x\in\Lreg$ then the map $L/\m_xL\to F/\m_xF$ is injective. Hence 
$L_x=\m_xL$ and $\dim L/\m_xL=c_L(X)$. Now (3) follows from (2) by 
globalization.  
\endproof

\proclaim
Theorem 5.2. The subset $\greg$ is open in $X$. Furthermore:

\item
(1) If $\greg=X$ then $A$ is a projective $A^{\g}$-module of rank 
$p^{c_\g(X)}$.

\item
(2) If $x\in\greg$ and $\n=\m_x\cap A^\g$ then $\n A$ is a maximal 
$\g$-invariant ideal of $A$ and the algebra $A/\n A$ is $\g$-equivariantly 
isomorphic with $\Hom_{u(\g_x)}\bigl(u(\g),k\bigr)$.

\item
(3) If $\greg=X$ then the assignment $I\mapsto I^\g$ establishes a bijection 
between the $\g$-invariant ideals of $A$ and the ideals of $A^\g$. 
%The inverse correspondence is given by $J\mapsto JA$. 
The canonical maps $A^\g\to(A/I)^\g$ are surjective.

\item
(4) If $\greg=X$ then $(B\otimes_{A^\g}A)^\g\cong B$ for every 
$A^\g$-algebra $B$ on which $\g$ operates trivially.

\item
(5) $[k(X):k(X)^{\g}]=p^{c_\g(X)}$.

\endproclaim

\Proof.
As is immediate from the definition of $L$, the $L$-invariants coincide with 
the $\g$-invariants, and an ideal of $A$ is stable under $L$ if and only 
if it is stable under $\g$. Since $\g_x=\rho^{-1}(L_x)$, the algebra map 
$\Hom_{u(L_x)}\bigl(u(L),\,k\bigr)\to\Hom_{u(\g_x)}\bigl(u(\g),\,k\bigr)$ 
induced by $\rho$ is an isomorphism. It follows that all assertions of the 
theorem for the $p$-Lie algebra $\g$ are equivalent to corresponding 
assertions for the $p$-Lie algebra $L$. Put $c=c_\g(X)=c_L(X)$.

Given $x\in\Lreg$, we have $\dim L_A^*/\m_xL_A^*=c$. Take  
$a_1,\ldots,a_c\in A$ such that $da_1,\ldots,da_c$ are a basis for 
a complement of $\m_xL_A^*$ in $L_A^*$. Since the $A_{\m_x}$-module 
$L_{\m_x}$ and its dual are free, passing to a suitable affine open 
neighbourhood of $x$, we may assume that $L$ is a free $A$-module and 
$da_1,\ldots,da_c$ are a basis for $L_A^*$ over $A$. Let $D_1,\ldots,D_c$ be 
the dual basis for $L$ over $A$. This means that 
$D_ia_l=da_l(D_i)=\delta_{il}$ for all $i,l$.  As $L$ is a Lie subalgebra, 
we have $[D_i,D_j]=\sum_{l=1}^cg_{ijl}D_l$ for certain $g_{ijl}\in A$. 
Applying the derivations on both sides of the equality to $a_l$, we deduce 
$g_{ijl}=0$. Since $L$ is closed under $p$-th powers, we have 
$D_i^p=\sum_{l=1}^ch_{il}D_l$ for certain $h_{il}\in A$. We deduce 
similarly that $h_{il}=0$. Thus the linear span $\a=\langle 
D_1,\ldots,D_c\rangle\subset L$ is an abelian Lie subalgebra with zero 
$p$-map. Since $L=A\a$, the assertions of the theorem for $L$ are equivalent 
to those for $\a$. Since $\dim\a<\infty$ they are equivalent also to the 
assertions of Theorem 2.1 for the corresponding action of the finite group 
scheme $\G(\a)$ (in fact this action is free).  
\endproof

\proclaim
Corollary 5.3.
Suppose that $\h\subset\g$ is a $p$-Lie subalgebra such that $\g=\h+\g_x$ 
for at least one $x\in\greg$. Then $A^\g=A^\h$.
\endproclaim

\Proof.
Obviously $\h_x=\g_x\cap\h$. By the hypotheses $\h/\h_x\cong\g/\g_x$ for 
some point $x\in\greg$. Then $c_\h(X)\ge\codim_\h\h_x=c_\g(X)$. It 
follows that $[K:K^\h]\ge[K:K^\g]$. On the other hand, $K^\g\subset K^\h$, 
whence $K^\g=K^\h$. We conclude that $A^\h=A\cap K^\h=A\cap K^\g=A^\g$.
\endproof

For every $r\ge1$ put $A^{(p^r)}=\{f^{p^r}\mid f\in A\}$. The notations 
$K^{(p^r)}$, $\m_x^{(p^r)}$ will have a similar meaning. For $f\in A$ let 
$d_xf:T_xX\to k$ denote the differential of $f$ at $x$.

\proclaim
Theorem 5.4.
Suppose that $X$ is a smooth affine variety and $f_1,\ldots,f_n\in A^\g$ 
where $n=\dim X-c_\g(X)$. Denote by $U$ the open set of those $x\in X$ for 
which $d_xf_1,\ldots,d_xf_n$ are linearly independent. If $U\ne\nothing$ 
then $K^\g=K^{(p)}(f_1,\ldots,f_n)$. If, moreover, $\codim_XX\setminus 
U\ge2$ then:

\item(1)
$A^\g=A^{(p)}[f_1,\ldots,f_n]$ and $A^\g$ is free of rank $p^n$ over 
$A^{(p)}$.

\item(2)
$A^\g$ is a locally complete intersection.

\item(3)
If $\pi:X\to X/\g=\Spec A^\g$ is the canonical morphism then $\pi(U)$ is the 
set of all smooth rational points of $X/\g$.

\endproclaim

\Proof.
Put $B=A^{(p)}[f_1,\ldots,f_n]$, $\,Y=\Spec B$, $\,X^{(p)}=\Spec A^{(p)}$.
The scheme $X^{(p)}$ is obtained from $X$ by base change $\f:k\to k$ where 
$\f$ is the Frobenius automorphism of $k$. Since smoothness is preserved 
under base change, $X^{(p)}$ is smooth. Denote by $\psi:Y\to X^{(p)}$ and 
$\phi:X\to Y$ the morphisms corresponding to the inclusions 
$A^{(p)}\subset B\subset A$. Both $\phi$ and $\psi$ are homeomorphisms. In 
particular, $X$, $Y$ and $X^{(p)}$ have the same dimension. Put $d=\dim X$.

For each commutative algebra $R$ denote by $\Omega_R$ the $R$-module of 
K\"ahler differentials of $R$ over $k$. By \cite{Ma}, (27.B),
$\dim_K\Omega_K=\degtr K/k=d$ since $K$ is separably generated over $k$. 
Furthermore, if $u_1,\ldots,u_d\in K$ are any elements such that 
$du_1,\ldots,du_d$ are a basis for $\Omega_K$ over $K$ then the elements 
$u_1^{m_1}\cdots u_d^{m_d}$ with $0\le m_i<p$ constitute a basis for $K$ 
over $K^{(p)}$. In particular, $[K:K^{(p)}]=p^d$. We have 
$\Omega_K\cong\Omega_A\otimes_AK$ where $\Omega_A$ is a projective 
$A$-module since $X$ is smooth (see \cite{Ma}, (29.B), Lemma 1). Assume that 
$U\ne\nothing$. If $x\in U$ then $df_1,\ldots,df_n$ are linearly 
independent modulo $\m_x\Omega_A$, hence constitute a basis for a direct 
summand of the free $A_{\m_x}$-module $(\Omega_A)_{\m_x}$. In particular, 
$df_1,\ldots,df_n$ are linearly independent over $A$, hence also over $K$. 
It follows that $[L:K^{(p)}]=p^n$ where $L=K^{(p)}(f_1,\ldots,f_n)$ is the 
field of fractions of $B$. Since $L\subset K^\g$ and 
$[K:L]=p^{d-n}=p^{c_\g(X)}=[K:K^\g]$, we deduce $L=K^\g$.

By the above $B$ is free over its subalgebra $A^{(p)}$ with basis elements 
$f_1^{m_1}\cdots f_n^{m_n}$ where $0\le m_i<p$. Then $B\cong 
A^{(p)}[t_1,\ldots,t_n]/I$ where $t_1,\ldots,t_n$ are indeterminates and $I$ 
is the ideal of the polynomial algebra generated by $n$ elements 
$t_i^p-f_i^p$. This means that $Y$ is isomorphic with the scheme-theoretic 
fibre $F=\tau^{-1}(0)$ of the morphism $\tau:X^{(p)}\times\An\to\An$ where 
$\An$ is the affine space of dimension $n$ and the components of $\tau$ are 
the functions $f_i^p-t_i^p$, $\,i=1,\ldots,n$. Note that $F$ is a complete 
intersection in the smooth variety $X^{(p)}\times\An$ since $F$ has
codimension $n$. It follows that $B$ is a locally complete intersection 
ring. In particular, $B$ is Cohen-Macaulay.

The tangent space $T_zF$ at a point $z\in F$ coincides with the kernel of 
the linear map $d_z\tau:T_z(X^{(p)}\times\An)\to T_0\An$ induced by $\tau$ 
in tangent spaces. Let $z=(x^{(p)},a)$ where $x^{(p)}=(\psi\circ\phi)(x)\in 
X^{(p)}$ for some $x\in X$ and $a\in\An$. It is easy to differentiate
$\tau$: for $(u,v)\in T_{x^{(p)}}X^{(p)}\oplus T_a\An\cong 
T_z(X^{(p)}\times\An)$ the vector $(d_z\tau)(u,v)\in T_0\An\cong k^n$ has 
components $(d_{x^{(p)}}f_i^p)(u)$, $\,i=1,\ldots,n$. Note that 
$T_{x^{(p)}}X^{(p)}\cong T_xX\otimes_\f k$ and the maps 
$d_{x^{(p)}}f_i^p:T_xX\otimes_\f k\to k$ are just $d_xf_i\otimes\id_k$ in 
this realization. The variety $F$ is smooth at $z$ if and only if $\dim 
T_zF=d$, if and only if $d_z\tau$ is surjective, if and only if 
$d_{x^{(p)}}f_1^p,\ldots,d_{x^{(p)}}f_n^p$ are linearly independent, if and 
only if $d_xf_1,\ldots,d_xf_n$ are linearly independent. In other words, 
the smoothness of $Y$ at $\phi(x)$ is equivalent to the inclusion $x\in U$. 
The codimension of the closed subset $Y\setminus\phi(U)$ in $Y$ is equal to 
that of $X\setminus U$ in $X$. Suppose that it is at least $2$. Then $Y$ is 
smooth in codimension $1$. By Serre's normality criterion $B$ is integrally 
closed (see \cite{Ma}, (17.I)). Then $A^\g=B$ since both algebras have the 
same field of fractions.  
\endproof

Suppose that $\g=\Lie\Gcal$ where $\Gcal$ is a reduced algebraic group 
operating on $X$ from the right. Then there is the induced action of $\g$ on 
$X$. For $x\in X$ denote by $\Gcal_x$ the scheme-theoretic stabilizer of $x$ 
in $\Gcal$. Let $\Greg\subset X$ be the open subset consisting of points 
$x$ for which the orbit $x\Gcal$ has a maximal possible dimension.

\proclaim
Theorem 5.5.
Suppose that $X$ is a smooth affine variety and $f_1,\ldots,f_n\in A^\Gcal$ 
where $n=\dim X-c_\g(X)$. Suppose also that the open set $U$ 
introduced in theorem {\rm5.4} is nonempty. Then 
$\greg=\{x\in\Greg\mid\Gcal_x$ is reduced$\}$. If $G$ denotes the $r$-th 
Frobenius kernel of $\Gcal$ for some $r\ge1$ then 
$K^G=K^{(p^r)}(f_1,\ldots,f_n)$. If, moreover, $\codim_XX\setminus U\ge2$ 
then:

\item(1)
$A^G=A^{(p^r)}[f_1,\ldots,f_n]$ and $A^G$ is free of rank $p^{rn}$ over 
$A^{(p^r)}$.

\item(2)
$A^G$ is a locally complete intersection.

\item(3)
$\pi(U)$ is the set of all smooth rational points of $X/G$ where $\pi:X\to 
X/G$ is the canonical morphism.

\endproclaim

\Proof.
Consider the morphism $\phi:X\to\An$ with components $f_1,\ldots,f_n$. Let 
$x\in U$. The differential $d_x\phi:T_xX\to k^n$ is then surjective.
By \cite{DG}, Ch.~I, \S4, 4.15 $\phi$ is smooth, hence also flat at $x$. 
Then $\dim_xF_x=\dim X-\dim\An=d-n$ by \cite{DG}, Ch.~I, \S3, 6.3, where 
$F_x=\phi^{-1}\bigl(\phi(x)\bigr)$ and $d=\dim X$. Since $F_x$ is a 
$\Gcal$-invariant closed subscheme of $X$, we have $x\Gcal\subset F_x$, and 
so $\dim x\Gcal\le d-n=c_\g(X)$. 

Since $U$ is open in $X$, we get $\dim x\Gcal\le c_\g(X)$, hence also 
$\dim\Gcal_x\ge\dim\Gcal-c_\g(X)$ for all $x\in X$. As $\dim\Gcal=\dim\g$, 
we can rewrite the last inequality in the form 
$$ 
c_\g(X)-\codim_\g\g_x\ge\dim\g_x-\dim\Gcal_x.\eqno(*)
$$
The subset $\Greg$ consists of those $x$ for which $\dim x\Gcal=c_\g(X)$, which 
is equivalent to an equality in $(*)$. By \cite{DG}, Ch.~III, \S2, 2.6 
$\Lie\Gcal_x=\g_x$. Furthermore, $\dim\g_x\ge\dim\Gcal_x$ and the equality 
holds here precisely when $\Gcal_x$ is smooth (which is equivalent to 
$\Gcal_x$ being reduced for an algebraically closed field) by \cite{DG}, 
Ch.~II, \S5, 2.1. As is easy to see, we have equalities everywhere above if 
and only if the left hand side of $(*)$ is zero, i.e., $x\in\greg$. This 
proves the first assertion of the theorem.

We have ${}|G|=p^{r\dim\g}$ (see \cite{Ja}, Part~I, 9.6, (2)). Next, 
$G_x=G\cap\Gcal_x$ for all $x\in X$. By \cite{Ja}, Part~I, 9.4, (2) $G_x$ 
coincides with the $r$-th Frobenius kernel of $\Gcal_x$. If $x\in\greg$ then 
$\Gcal_x$ is reduced, whence $|G_x|=p^{r\dim\g_x}$, and 
$(G:G_x)=|G|/|G_x|=p^{r\codim_\g\g_x}$. Since both $\Xreg$ and $\greg$ are 
open and nonempty, they have a common point, which shows that 
$q(X)=p^{rc_\g(X)}$ in the notations of Theorem 2.1. We deduce that
$[K:K^G]=p^{rc_\g(X)}$.

Let $L_i=K^{(p^i)}(f_1,\ldots,f_n)$ for each $i\ge1$. If $i>1$ then
$K^{(p^i)}\subset L_{i-1}^{(p)}\subset L_i$. As we have seen in Theorem 5.4 
the elements $f_1^{m_1}\cdots f_n^{m_n}$ with $0\le m_i<p$ are linearly 
independent over $K^{(p)}$. It follows that $[L_i:L_{i-1}^{(p)}]=p^n$. Hence
$$
[L_i:K^{(p^i)}]=
p^n[L_{i-1}^{(p)}:K^{(p^i)}]=p^n[L_{i-1}:K^{(p^{i-1})}].
$$ 
We have also $[K:K^{(p^i)}]=[K:K^{(p)}]\cdot[K^{(p)}:K^{(p^i)}]
=p^d[K:K^{(p^{i-1})}]$. Proceeding by induction on $i$ we 
deduce that $[L_i:K^{(p^i)}]=p^{ni}$ and $[K:K^{(p^i)}]=p^{di}$. Taking 
$i=r$, we get $[K:L_r]=p^{(d-n)r}=[K:K^G]$. Since $L_r\subset K^G$, it 
follows $L_r=K^G$. The remainder of the theorem is proved in the same way as 
Theorem 5.4 with obvious changes. 
\endproof

Two classical cases of Theorem 5.5 are those when $\Gcal$ is a semisimple 
algebraic group operating either on $\g$ via the adjoint representation or 
on itself by conjugations. Let $n$ denote the rank of $\Gcal$. Under 
assumption that $p$ does not divide the order of the Weyl group of $\Gcal$ 
it was shown by Veldkamp \cite{Veld} that $k[\g]^\Gcal$ is generated by $n$ 
algebraically independent polynomials $J_1,\ldots,J_n$ and $\g\reg$ consists 
precisely of those $x\in\g$ for which $d_xJ_1,\ldots,d_xJ_n$ are linearly 
independent. The complement of $\g\reg$ has pure codimension $3$ in $\g$.
The stabilizer $\g_x$ is just the centralizer of $x$ in $\g$.
As $\dim\g_x=n$ for all $x\in\g\reg$, we have $c_\g(\g)=\dim\g-n$.

Suppose that $\Gcal$ is simply connected. Then the algebra of regular 
functions on $\Gcal$, constant on the conjugacy classes, is generated by 
the characters $\chi_1,\ldots,\chi_n$ of fundamental irreducible 
representations of $\Gcal$. As shown by Steinberg \cite{St}, $\Gcal\reg$ 
consists precisely of those $x\in\Gcal$ for which 
$d_x\chi_1,\ldots,d_x\chi_n$ are linearly independent. The complement of 
$\Gcal\reg$ in $\Gcal$ again has codimension 3. If $x\in\Gcal$ then 
$\g_x=\ker d_e\mu_x$ where $\mu_x:\Gcal\to\Gcal$ is the morphism defined by 
the rule $y\mapsto y^{-1}xy$. It follows $\g_x=\{v\in\g\mid(\Ad x)v=v\}$. If 
$\T\subset\Gcal$ is a maximal torus then there exists $t\in\T$ such that 
$\alpha(t)\ne1$ for all roots $\alpha$. Then $\g_t=\Lie\T$, and so 
$c_\g(\Gcal)=\dim\Gcal-n$. No restrictions on $p$ are needed in this case.  
The invariants of $\g$ and Frobenius kernels were described by Friedlander 
and Parshall \cite{Fri}, and Donkin \cite{Don}.

We give yet another example showing that Theorem 5.4 has a wider range of 
applications. Let $\g=W_n$ be the Jacobson-Witt algebra. Recall that 
$\g=\Der B_n$ where $B_n=k[x_1,\ldots,x_n]$, $\,x_i^p=0$, is the truncated 
polynomial algebra. If $\Gcal=\nobreak\Aut B_n$ then $\Lie\Gcal$ is the 
subalgebra of codimension $n$ in $\g$ consisting of derivations that leave 
stable the maximal ideal $\n$ of $B_n$. For $D\in\g$ denote by $\chi_D(t)$ 
the characteristic polynomial of $D$ as a linear transformation of $B_n$.  
As is proved by Premet in \cite{Pre}, 
$\chi_D(t)=t^{p^n}+\sum_{i=0}^{n-1}\psi_i(D)t^{p^i}$ where $\psi_i$ are 
algebraically independent polynomial functions generating the algebra 
$k[\g]^\Gcal$. There exists an open subset of $\g$ consisting of elements 
$D$ such that $\dim\g_D=n$ and $\g=\Lie\Gcal\oplus\g_D$ where $\g_D$ 
is the centralizer of $D$ in $\g$ (\cite{Pre}, Lemma 1). Hence
$c_\g(\g)=d-n$ where $d=\dim\g$, and $\psi_0,\ldots,\psi_{n-1}$ are 
$\g$-invariant according to Corollary 5.3. Let $\phi:\g\to\An$ be the 
morphism with components $\psi_0,\ldots,\psi_{n-1}$ and $U$ the subset of 
those $D\in\g$ for which $d_D\psi_0,\ldots,d_D\psi_{n-1}$ are linearly 
independent. Premet proved that each fibre 
$F_D=\phi^{-1}\bigl(\phi(D)\bigr)$ is irreducible of dimension $d-n$ and 
$F_D\cap U$ is nonempty (\cite{Pre}, Lemmas 12 and 13). 

Put $U_1=\{D\in\g\mid\psi_0(D)\ne0\}$. If $D\in U_1$, then $D$ is a linear 
combination of $D^p,\ldots,D^{p^n}$ as 
$D^{p^n}+\sum_{i=0}^{n-1}\psi_i(D)D^{p^i}=0$. Hence $D$ is semisimple, and 
we can find its eigenvectors $y_1,\ldots,y_n\in B_n$ such that 
$B_n=\langle y_1,\ldots,y_n\rangle\oplus(k+\n^2)$. Then the monomials 
$y_1^{m_1}\cdots y_n^{m_n}$ with $0\le m_i<p$ constitute a basis for $B_n$.  
Let $\lambda_i$ be the eigenvalue of $y_i$. Since the rank of $D$ as a 
linear transformation is equal to $p^n-1$, the equality 
$\sum m_i\lambda_i=0$ can hold only for $m_1=\ldots=m_n=0$. This implies that 
$D$ generates a torus of dimension $n$ in $\g$. If $D'\in\g_D$ then 
$$ 
\textstyle\sum\limits_{i=0}^{n-1}(d_D\psi_i)(D')D^{p^i}=-\psi_0(D)D',
$$
which is a special case of \cite{Pre}, Lemma 7, (i). Taking $D'=D^{p^j}$ 
where $0\le j<n$, we see that $(d_D\psi_i)(D')\ne0$ only for $i=j$. Hence 
$U_1\subset U$. Suppose that $Z$ is an irreducible component of the closed 
subset $\g\setminus U_1$ having codimension 1 in $\g$ (in fact it can be 
shown that $\psi_0$ is irreducible). Note that $\overline{\phi(Z)}\ne\An$ as 
$\psi_0(Z)=\{0\}$. By the theorem on dimensions of fibres we have $\dim 
Z\cap F_D\ge\dim Z-\dim\overline{\phi(Z)}\ge d-n$ for all $D\in Z$. It 
follows that $Z$ is a union of fibres of $\phi$. In particular, $Z\cap 
U\ne\nothing$. We deduce that $\codim_\g\g\setminus U\ge2$. Thus we meet the 
hypotheses of Theorem 5.4:

\proclaim
Corollary 5.6.
If $A$ denotes the algebra of polynomial functions on $W_n$ then 
$A^{W_n}=A^{(p)}[\psi_0,\ldots,\psi_{n-1}]$. Moreover, $A^{W_n}$ is 
free of rank $p^n$ over $A^{(p)}$ and is a locally complete intersection.  
\endproclaim

In conclusion we make comments concerning the results of section 4. Assume 
that $\dim\g<\infty$. According to \cite{Hoch} the algebra $u(\g_x)$ is 
semisimple if and only if $\g_x$ is a torus. Thus Theorem 4.2 says that 
$\g_x$ is toral if and only if $x$ lies in an affine open subset $U\subset 
X$ such that $k[U]$ is an injective $u(\g)$-module. Such points $x$ exist
if and only if $k(X)$ is an injective $u(\g)$-module. A.~Premet pointed out 
to me that the openness of the set of points with a toral $\g_x$ can be 
proved by geometric arguments.

\references
\tracingmacros=1
\nextref
Bou
\auth
N. Bourbaki     
\endauth
\book{Commutative Algebra}
\publisher{Springer, Berlin}
\Year{1989}

\nextref
Cline
\auth
E. Cline B. Parshall L. Scott 
\endauth
\paper{A Mackey imprimitivity theory for algebraic groups}
\journal{Math.~Z.}
\Vol{182}
\Year{1983}
\Pages{447--471}

\nextref
DG
\auth
M. Demazure P. Gabriel   
\endauth
\book{Groupes Alg\'ebriques I}
\publisher{Masson, Paris}
\Year{1970}

\nextref
Doi
\auth
Y. Doi     
\endauth
\paper{Algebras with total integrals}
\journal{Comm. Algebra}
\Vol{13}
\Year{1985}
\Pages{2137--2159}

\nextref
Don
\auth
S. Donkin     
\endauth
\paper{Infinitesimal invariants of algebraic groups}
\journal{J.~London Math. Soc.}
\Vol{45}
\Year{1992}
\Pages{481--490}

\nextref
Fri
\auth
E.M. Friedlander B.J. Parshall   
\endauth
\paper{Rational actions associated to the adjoint representation}
\journal{Ann. Sci. \'Ecole Norm. Sup.}
\Vol{20}
\Year{1987}
\Pages{215--226}

\nextref
Hoch
\auth
G.P. Hochschild     
\endauth
\paper{Representations of restricted Lie algebras of characteristic $p$}
\journal{Proc. Amer. Math. Soc.}
\Vol{5}
\Year{1954}
\Pages{603--605}

\nextref
Hoch2
\auth
G.P. Hochschild     
\endauth
\paper{Simple algebras with purely inseparable splitting fields of exponent 1}
\journal{Trans. Amer. Math. Soc.}
\Vol{79}
\Year{1955}
\Pages{477--489}

\nextref
Ja
\auth
J. Jantzen     
\endauth
\book{Representations of algeraic groups}
\publisher{Academic Press}
\Year{1987}

\nextref
Kop
\auth
M. Koppinen T. Neuvonen   
\endauth
\paper{An imprimitivity theorem for Hopf algebras}
\journal{Math. Scand.}
\Vol{41}
\Year{1977}
\Pages{193--198}

\nextref
Ma
\auth
H. Matsumura     
\endauth
\book{Commutative Algebra, Second Edition}
\publisher{Benjamin}
\Year{1980}

\nextref
Mil
\auth
A.A. Mil'ner     
\endauth
\paper{Irreducible representations of modular Lie algebras}
\journal{Izv. Akad. Nauk SSSR Ser. Mat.}
\Vol{39}
\Year{1975}
\Pages{1240--1259}
[In Russian]
Translation in \journal{Math. USSR Izv.}
\Vol{9}
\Year{1975}
\Pages{1169--1187}

\nextref
Mon
\auth
S. Montgomery     
\endauth
\book{Hopf algebras and Their Actions on Rings}
\bookseries{CBMS Regional Conference Series in Mathematics}
\Vol{82}
\publisher{American Mathematical Society}
\Year{1993}

\nextref
Mum
\auth
D. Mumford     
\endauth
\book{Abelian Varieties}
\publisher{Oxford University Press, Oxford}
\Year{1970}

\nextref
Nich
\auth
W.D. Nichols M.B. Zoeller   
\endauth
\paper{A Hopf algebra freeness theorem}
\journal{Amer. J. Math.}
\Vol{111}
\Year{1989}
\Pages{381--385}

\nextref
Ober
\auth
U. Oberst H.-J. Schneider   
\endauth
\paper{\"Uber Untergruppen endlicher algebraischer Gruppen}
\journal{Manuscripta Math.}
\Vol{8}
\Year{1973}
\Pages{217--241}

\nextref
Pre
\auth
A.A. Premet     
\endauth
\paper{The theorem on restriction of invariants and nilpotent elements in $W_n$}
\journal{Mat. Sbornik}
\Vol{182}
\Year{1991}
\Pages{746--773}
[In Russian]
Translation in \journal{Math. USSR Sbornik}
\Vol{73}
\Year{1992}
\Pages{135--159}

\nextref
St
\auth
R. Steinberg     
\endauth
\paper{Regular elements of semisimple algebraic groups}
\journal{IHES Publ. Math.}
\Vol{25}
\Year{1965}
\Pages{49--80}

\nextref
Sw
\auth
M.E. Sweedler     
\endauth
\book{Hopf Algebras}
\publisher{Benjamin, New York}
\Year{1969}

\nextref
Tak
\auth
M. Takeuchi     
\endauth
\paper{Relative Hopf modules---equivalences and freeness criteria}
\journal{J. Algebra}
\Vol{60}
\Year{1979}
\Pages{452--471}

\nextref
Veld
\auth
F.D. Veldkamp     
\endauth
\paper{The center of the universal enveloping algebra of a Lie algebra in characteristic $p$}
\journal{Ann. Sci. \'Ecole Norm. Sup.}
\Vol{5}
\Year{1972}
\Pages{217--240}

\nextref
Wat
\auth
W.C. Waterhouse     
\endauth
\paper{Geometrically reductive affine group schemes}
\journal{Arch. Math.}
\Vol{62}
\Year{1994}
\Pages{306--307}

\endreferences
\bye